\documentclass[a4paper,12pt]{amsart}
\usepackage{amsfonts}
\usepackage{amssymb}
\usepackage{ifthen}
\usepackage{graphicx}
\nonstopmode \numberwithin{equation}{section}
\usepackage[usenames]{color}
\newtheorem{thm}{Theorem}
\newtheorem{lem}{Lemma}
\newtheorem{cor}{Corollary}[section]

\newtheorem{cl}{Claim}
\newtheorem{ca}{Case}
\newtheorem{sca}{Subcase}
\newtheorem{scl}{Subclaim}
\newtheorem{conj}[equation]{Conjecture}

\theoremstyle{definition}
\newtheorem{defn}{Definition}

\newtheorem{op}[equation]{Open Problem}
\newtheorem{ques}[equation]{Question}
\newtheorem{rem}{Remark}[section]
\newtheorem{exam}[equation]{Example}

\newcounter {own}
\def\theown {\thesection       .\arabic{own}}

\newenvironment{pf}[1][]{%
 \vskip 3mm
 \noindent
 \ifthenelse{\equal{#1}{}}%
  {{\slshape Proof. }}%
  {{\slshape #1.} }%
 }%
{\qed\bigskip}

\newcounter{alphabet}
\newcounter{tmp}
\newenvironment{Thm}[1][]{\refstepcounter{alphabet}%
\bigskip%
\noindent%
{\bf Theorem \Alph{alphabet}}%
\ifthenelse{\equal{#1}{}}{}{ (#1)}%
{\bf .} \itshape}{\vskip 8pt}

\makeatletter
\newcommand{\Ref}[1]{\@ifundefined{r@#1}{}{\setcounter{tmp}{\ref{#1}}\Alph{tmp}}}
\makeatother

\def\be{\begin{equation}}
\def\ee{\end{equation}}

\newcommand{\bee}{\begin{enumerate}}
\newcommand{\eee}{\end{enumerate}}

\newcommand{\blem}{\begin{lem}}
\newcommand{\elem}{\end{lem}}
\newcommand{\bthm}{\begin{thm}}
\newcommand{\ethm}{\end{thm}}
\newcommand{\bcor}{\begin{cor}}
\newcommand{\ecor}{\end{cor}}
\newcommand{\beg}{\begin{exam}}
\newcommand{\eeg}{\end{exam}}
\newcommand{\begs}{\begin{examples}}
\newcommand{\eegs}{\end{examples}}
\newcommand{\bdefe}{\begin{defn}}
\newcommand{\edefe}{\end{defn}}
\newcommand{\bprob}{\begin{prob}}
\newcommand{\eprob}{\end{prob}}
\newcommand{\bques}{\begin{ques}}
\newcommand{\eques}{\end{ques}}
\newcommand{\bei}{\begin{itemize}}
\newcommand{\eei}{\end{itemize}}
\newcommand{\bcon}{\begin{conj}}
\newcommand{\econ}{\end{conj}}
\newcommand{\bop}{\begin{op}}
\newcommand{\eop}{\end{op}}

\newcommand{\bca}{\begin{ca}}
\newcommand{\eca}{\end{ca}}
\newcommand{\bsca}{\begin{sca}}
\newcommand{\esca}{\end{sca}}

\newcommand{\bcl}{\begin{cl}}
\newcommand{\ecl}{\end{cl}}

\newcommand{\bscl}{\begin{scl}}
\newcommand{\escl}{\end{scl}}

\newcommand{\bcons}{\begin{conjs}}
\newcommand{\econs}{\end{conjs}}
\newcommand{\bprop}{\begin{propo}}
\newcommand{\eprop}{\end{propo}}
\newcommand{\br}{\begin{rem}}
\newcommand{\er}{\end{rem}}
\newcommand{\brs}{\begin{rems}}
\newcommand{\ers}{\end{rems}}
\newcommand{\bo}{\begin{obser}}
\newcommand{\eo}{\end{obser}}
\newcommand{\bos}{\begin{obsers}}
\newcommand{\eos}{\end{obsers}}
\newcommand{\bpf}{\begin{pf}}
\newcommand{\epf}{\end{pf}}
\newcommand{\ba}{\begin{array}}
\newcommand{\ea}{\end{array}}
\newcommand{\beq}{\begin{eqnarray}}
\newcommand{\beqq}{\begin{eqnarray*}}
\newcommand{\eeq}{\end{eqnarray}}
\newcommand{\eeqq}{\end{eqnarray*}}

\newcounter{minutes}
\divide\time by 60
\newcounter{hours}
\multiply\time by 60 \addtocounter{minutes}{-\time}

\begin{document}
\bibliographystyle{amsplain}
\title [] {Characterizations of Hardy-type, Bergman-type and Dirichlet-type
spaces on certain classes of complex-valued functions}

\def\thefootnote{}
\footnotetext{ \texttt{\tiny File:~\jobname .tex,
          printed: \number\day-\number\month-\number\year,
          \thehours.\ifnum\theminutes<10{0}\fi\theminutes}
} \makeatletter\def\thefootnote{\@arabic\c@footnote}\makeatother


\author{Shaolin Chen}
\address{Shaolin Chen, Department of Mathematics and Computational Science,
Hengyang Normal University, Hengyang, Hunan 421008, People's
Republic of China.} \email{mathechen@126.com}

\author{Antti Rasila }
\address{Antti Rasila, Department of Mathematics and Systems Analysis, Aalto University, P. O. Box 11100, FI-00076 Aalto,
 Finland.} \email{antti.rasila@iki.fi}

 \author{Matti Vuorinen}
\address{Matti Vuorinen,
Department of Mathematics, University of Turku, Turku 20014,
Finland. }\email{vuorinen@utu.fi}



\subjclass[2000]{Primary: 32A10, 30D55; Secondary: 30C65, 58J10}
\keywords{Hardy-type space, Bergman-type space, Dirichlet-type
space.}

\begin{abstract}
In this paper, we continue our investigation of function spaces  on
certain classes of complex-valued functions. In particular, we give
characterizations on Hardy-type, Bergman-type and Dirichlet-type
spaces. Furthermore, we present applications of our results to
certain nonlinear PDEs.
\end{abstract}

\maketitle \pagestyle{myheadings} \markboth{Sh. Chen,   A. Rasila
and M. Vuorinen }{ Characterizations of Hardy-type, Bergman-type and
Dirichlet-type spaces}

\section{Introduction and main results}\label{csw-sec1}
For a positive integer $n\geq1$, let $\mathbb{C}^{n}$ denote the
complex {\it Euclidean $n$-space}. For $z:=(z_{1},\ldots,z_{n})$ and
$w=(w_{1},\ldots,w_{n})$ in $\mathbb{C}^{n}$, we let
$\overline{z}=(\overline{z}_{1},\ldots,\overline{z}_{n} ),$ and
$\langle z,w\rangle := \sum_{k=1}^nz_k\overline{w}_k$ with the
 {\it Euclidean norm} $ \|z\|:={\langle z,z\rangle}^{1/2}$
which makes $\mathbb{C}^n$ into an $n$-dimensional complex {\it
Hilbert} space. For $a\in \mathbb{C}^n$ and $r>0$,
$\mathbb{B}^n(a,r)$
denotes the (open) ball of radius $r$ with center $a$. Also, we let
$\mathbb{B}^n(r):=\mathbb{B}^n(0,r)$ and
denote the unit ball by $\mathbb{B}^n:=\mathbb{B}^n(1)$. In particular, let
$\mathbb{B}^1(r)=\mathbb{D}(r)$ and  $\mathbb{D}=\mathbb{B}^1$. For
a domain $\Omega\subset\mathbb{C}^{n}$ with non-empty boundary, let
$d_{\Omega}(z)$ be the Euclidean distance from $z$ to the boundary
$\partial\Omega$ of $\Omega$. Moreover, we always use $d(z)$ to
denote the Euclidean distance from $z$ to the boundary of
$\mathbb{B}^{n}$. We denote by $\mathcal{C}^{m}(\mathbb{B}^{n})$ the set of all
$m$-time continuously differentiable complex-valued functions $f$ of
$\mathbb{B}^{n}$ into $\mathbb{C}$, where $m\in\{0,1,\ldots\}$.

For $k\in\{1,\ldots,n\},$ $z=(z_{1},\ldots,z_{n})\in\mathbb{C}^{n}$
and $f\in\mathcal{C}^{1}(\mathbb{B}^{n})$, we introduce the
following notations:
$$\nabla f=(f_{z_{1}},\ldots,f_{z_{n}}),~
\overline{\nabla}f=(f_{\overline{z}_{1}},\ldots,f_{\overline{z}_{n}})~\mbox{and}~D_{f}=(\nabla
f,~\overline{\nabla}f),
$$ where
$f_{z_{k}}=\partial f/\partial z_{k}=1/2\big(\partial f/\partial
x_{k}-i\partial f/\partial y_{k}\big)$,
$f_{\overline{z}_{k}}=\partial f/\partial
\overline{z}_{k}=1/2\big(\partial f/\partial x_{k}+i\partial
f/\partial y_{k}\big)$ and $z_{k}=x_{k}+iy_{k}$, with $x_{k}$ and
$y_{k}$ real.
 Let
 $\|D_{f}\|$ be the {\it Hilbert-Schmidt semi-norm } given by $$\|D_{f}\|=(\|\nabla
f\|^{2}+\|\overline{\nabla}f\|^{2})^{1/2}.$$

Let $f=u+iv\in\mathcal{C}^{1}(\mathbb{B}^{n})$, where $u$ and $v$
are real-valued functions. Then for
$z=(z_{1},\ldots,z_{n})=(x_{1}+iy_{1},\ldots,x_{n}+iy_{n})\in\mathbb{B}^{n}$,
\be\label{1x} \|\nabla f(z)\|+\|\overline{\nabla}f(z)\|\leq \|\nabla
u(z)\|+\|\nabla v(z)\|, \ee where $$\nabla u=\Big( \frac{\partial
u}{\partial x_{1}}, \frac{\partial u}{\partial y_{1}},\ldots,
\frac{\partial u}{\partial x_{n}}, \frac{\partial u}{\partial
y_{n}}\Big)~\text{ and }~ \nabla v=\Big( \frac{\partial v}{\partial
x_{1}}, \frac{\partial v}{\partial y_{1}},\ldots, \frac{\partial
v}{\partial x_{n}}, \frac{\partial v}{\partial y_{n}} \Big).$$ Note
that the converse of (\ref{1x}) is not always true (see
\cite{CPW3}).


\subsection*{Generalized Hardy spaces} For $p\in(0,\infty]$, the {\it generalized Hardy space}
$\mathcal{H}_{g}^{p}(\mathbb{B}^{n})$
 consists of measurable functions $f:\ \mathbb{B}^{n}\rightarrow\mathbb{C}$
 such that  $M_{p}(r,f)$ exists for all $r\in(0,1)$ and  $ \|f\|_{p}<\infty$, where
$$
\|f\|_{p}=
\begin{cases}
\displaystyle\sup_{0<r<1}M_{p}(r,f),
& \mbox{ if } p\in(0,\infty),\\
\displaystyle\sup_{z\in\mathbb{B}^{n}}|f(z)|, &\mbox{ if }\,
p=\infty,
\end{cases}~ M_{p}(r,f)=\left(\int_{\partial\mathbb{B}^{n}}|f(r\zeta)|^{p}\,d\sigma(\zeta)\right)^{1/p}
$$
and $d\sigma$ denotes the normalized Lebesgue surface measure in
$\partial\mathbb{B}^{n}$.

There are numerous characterizations of the classical analytic Hardy spaces in the literature, see for example \cite{Du1,GP1,GPP, HL1, HL2, Pav}. But, to our knowledge, there are few analogous results for general complex-valued
functions. In this paper, we give the following characterization of a class of
complex-valued functions $f$ in Hardy-type spaces.

\begin{thm}\label{thm-3}
For $p\geq2$, let $f\in\mathcal{C}^{2}(\mathbb{B}^{n})$ with ${\rm
Re}(f\overline{\Delta f})\geq0.$ Then,
$$\int_{\mathbb{B}^{n}}d(z)\Delta\big(|f(z)|^{p}\big)dV_{N}(z)<\infty$$
if and only if  $f\in \mathcal{H}_{g}^{p}(\mathbb{B}^{n})$, where
$\Delta$ is the usual complex Laplacian operator
$$\Delta:=4\sum_{k=1}^{n}\frac{\partial^{2}}{\partial z_{k}\partial \overline{z}_{k}}=\sum_{k=1}^{n}
\left(\frac{\partial^{2}}{\partial
x_{k}^{2}}+\frac{\partial^{2}}{\partial y_{k}^{2}}\right)$$ for
$z=(z_{1},\ldots,z_{n})=(x_{1}+iy_{1},\ldots,x_{n}+iy_{n})\in\mathbb{B}^{n}$.
\end{thm}

\subsection*{Yukawa PDE}
Let $\tau,~\eta:~\mathbb{B}^{n}\rightarrow[0,\infty)$ be continuous
and $f=u+iv\in\mathcal{C}^{2}(\mathbb{B}^{n})$, where $u$ and $v$
are real-valued functions in $\mathbb{B}^{n}$. The nonlinear
elliptic partial differential equation (PDE) of the
form \be\label{eq-CR1.1} \Delta f(z)=\tau(z)
f(z)+\eta(z)\mbox{Re}\big(f(z)\big) \ee is called the {\it
non-homogeneous Yukawa PDE}, where $z\in\mathbb{B}^{n}$.  If $\tau$
in (\ref{eq-CR1.1}) is a positive constant function and
$\eta\equiv0$, then we have the usual Yukawa PDE. This equation
arose from the work of the Japanese Nobel physicist Hideki Yukawa,
who used it to describe the nuclear potential of a point charge as
$e^{-\sqrt{\tau} r}/r$  (cf. \cite{A,CPR,CRW,Do,D-,D-1,E,SW,Ya}).

As an application of Theorem \ref{thm-3}, we obtain the following
result.

\begin{cor}\label{cor-1}
For $p\geq2$, let $f\in\mathcal{C}^{2}(\mathbb{B}^{n})$ satisfying
{\rm (\ref{eq-CR1.1})}. Then,
$$\int_{\mathbb{B}^{n}}d(z)\Delta\big(|f(z)|^{p}\big)dV_{N}(z)<\infty$$
if and only if  $f\in \mathcal{H}_{g}^{p}(\mathbb{B}^{n})$.
\end{cor}

A continuous increasing function $\omega:\, [0,\infty)\rightarrow
[0,\infty)$ with $\omega(0)=0$ is called a {\it majorant} if
$\omega(t)/t$ is non-increasing for $t>0$ (cf.
\cite{Dy1,Dy2,P,Pav1}). Given a subset $\Omega$ of $\mathbb{C}$, a
function $f:\, \Omega\rightarrow \mathbb{C}$ is said to belong to
the {\it Lipschitz space $L_{\omega}(\Omega)$} if there
is a positive constant $C$ such that 
$$|f(z)-f(w)|\leq C\omega(|z-w|) ~\mbox{ for all $z,w\in\Omega.$} $$ 

A classical result of Hardy and Littlewood asserts that if
$p\in(0,\infty]$, $\alpha\in(1,\infty)$ and $f$ is an analytic
function in $\mathbb{D}$, then (cf. \cite{Du1,HL1,HL2})
$$ M_{p}(r,f')=O \left(\Big(\frac{1}{1-r}\Big)^{\alpha} \right ) ~\mbox{ as $r\rightarrow1$},
$$
if and only if
$$M_{p}(r,f)=O \left (\Big(\log\frac{1}{1-r}\Big )^{\alpha-1}\right) ~\mbox{ as $r\rightarrow1$}.
$$

In \cite{GPP}, via the closed graph theorem, Girela, Pavlovi\'c and
Pel\'{a}ez refined the above result for the case $\alpha=1$ as
follows.

%
\begin{Thm}\label{Thm-0.1} $($\cite[Theorem 1.1]{GPP}$)$
\label{ThmB} Let $p\in(2,\infty)$. For $r\in(0,1)$, if  $f$ is
analytic in $\mathbb{D}$ such that
$$M_{p}(r,f')=O \left ( \frac{1}{1-r} \right  )  ~\mbox{ as $r\rightarrow1$},
$$
then
$$ M_{p}(r,f)=O \left (\Big(\log\frac{1}{1-r}\Big)^{\frac{1}{2}} \right ) ~\mbox{ as $r\rightarrow 1$}
$$ and the exponent $1/2$ is sharp.
\end{Thm}

Theorem \Ref{Thm-0.1} gives an affirmative  answer to the open
problem in \cite[p. 464, Equation (26)]{GP1}. For related
investigations on this topic, we refer to \cite{CPR,CPW2,CRW}.

Next we study the relationship between the
 integral means of solutions to the
equation (\ref{eq-CR1.1}) and those of their two order partial
derivative. Our result is given as follows.

\begin{thm}\label{thm-5.0}
Let $\omega$ be a majorant
 and $f\in\mathcal{C}^{2}(\mathbb{B}^{n})$ satisfying {\rm
(\ref{eq-CR1.1})} with $\eta+\tau<4n/p$, where $\tau$ and $\eta$ are
nonnegative constant functions. For $p\geq2$ and $r\in(0,1)$, if
$$M_{p}(r, D^{\ast}_{f})\leq
M^{\ast}\omega\left(\frac{1}{1-r}\right),$$ then
$$M_{p}(r,D_{f})\leq
\sqrt{M^{\ast}_{2}}
\left[\|D_{f}(0)\|^{2}+M^{\ast}_{1}\int_{0}^{1}\omega\Big(\frac{1}{1-rt}\Big)dt\right]^{\frac{1}{2}},$$
and   $f\in \mathcal{H}_{g}^{p}(\mathbb{B}^{n})$, where $M^{\ast}$
is a positive constant,
$$D^{\ast}_{f}=\left[\sum_{j=1}^{n}\sum_{k=1}^{n}\Big(|f_{z_{k}z_{j}}|^{2}+|f_{z_{k}\overline{z}_{j}}|^{2}+
|f_{\overline{z}_{k}z_{j}}|^{2}+
|f_{\overline{z}_{k}\overline{z}_{j}}|^{2}\Big)\right]^{\frac{1}{2}},$$
$M^{\ast}_{1}=2p(2p-3)(M^{\ast})^{2}\omega(1)$ and
$M^{\ast}_{2}=1/\left[1-p(\eta+\tau)/(4n)\right].$
\end{thm}

In particular, by taking $\omega(t)=t$ in Theorem \ref{thm-5.0}, we
obtain the following result.

\begin{cor}\label{cor-5.0}
Let $p\geq2$ and $f\in\mathcal{C}^{2}(\mathbb{B}^{n})$ satisfying
{\rm (\ref{eq-CR1.1})} with $\eta+\tau<4n/p$, where $\tau$ and
$\eta$ are nonnegative constant functions. For $r\in(0,1)$, if

$$M_{p}(r, D^{\ast}_{f})=O \left ( \frac{1}{1-r} \right  )  ~\mbox{ as $r\rightarrow1$},
$$
then
$$ M_{p}(r,D_{f})=O \left (\Big(\log\frac{1}{1-r}\Big)^{\frac{1}{2}} \right ) ~\mbox{ as $r\rightarrow
1$},
$$
and  $f\in \mathcal{H}_{g}^{p}(\mathbb{B}^{n})$.
\end{cor}

\subsection*{Dirichlet-type spaces
and Bergman-type spaces} For $\nu,~\mu,~t\in\mathbb{R}$,
$$\mathcal{D}_{f}(\nu,\mu,t)=\int_{\mathbb{B}^{n}}d^{\nu}(z)|f(z)|^{\mu}\|D_{f}(z)\|^{t}dV_{N}(z)<\infty$$
is called {\it Dirichlet-type energy integral} of
 the complex-valued function $f$,
where $dV_{N}$ denotes the normalized Lebesgue volume measure in
$\mathbb{B}^{n}$ (cf. \cite{AIM,A,CRW,E,GPP,GP,SH,ST,W,Ya}). In
particular, for $\nu\geq0$, $\mu=0$ and $0<t<\infty$, we use
$\mathcal{D}_{\nu,t}(\mathbb{B}^{n})$ to denote the {\it
Dirichlet-type space} consisting of  all
$f\in\mathcal{C}^{1}(\mathbb{B}^{n})$ with the norm
$$\|f\|_{\mathcal{D}_{\nu,t}}=|f(0)|+\big(\mathcal{D}_{f}(\nu,0,t)\big)^{1/t}<\infty.
$$
Moreover, for $\nu>-1$, $0<\mu<\infty$ and $t=0$, we denote by
$b_{\nu,\mu}(\mathbb{B}^{n})$ the {\it Bergman-type space}
consisting of all $f\in\mathcal{C}^{0}(\mathbb{B}^{n})$ with the
norm
$$\|f\|_{b_{\nu,\mu}}=|f(0)|+\big(\mathcal{D}_{f}(\nu,\mu,0)\big)^{1/\mu}<\infty.$$

We refer to \cite{DS,GPP,GP,GP-1,Pav-1,Pav1,Z} for basic
characterizations of analytic (or harmonic) Bergman-type spaces and
Dirichlet-type spaces.  Again, for general complex-valued functions, very little related
research can be found from the literature. The following is  a
characterization of  a class of complex-valued functions $f$ in
Bergman-type spaces.

\begin{thm}\label{thm-2}
 Let $f\in\mathcal{C}^{2}(\mathbb{B}^{n})$ with ${\rm
Re}(f\overline{\Delta f})\geq0.$ Then,  for $p\geq2$ and
$\alpha\geq2$,
$$\int_{\mathbb{B}^{n}}(1-|z|^{2})^{\alpha}\Delta\big(|f(z)|^{p}\big)dV_{N}(z)<\infty,$$
if and only if  $f\in b_{\alpha-2,p}(\mathbb{B}^{n})$.
\end{thm}

The following result easily follows from Theorem \ref{thm-2}.

\begin{cor}\label{cor-2}
Let $f\in\mathcal{C}^{2}(\mathbb{B}^{n})$ satisfying {\rm
(\ref{eq-CR1.1})}. Then, for $p\geq2$ and  $\alpha\geq2$,
$$\int_{\mathbb{B}^{n}}(1-|z|^{2})^{\alpha}\Delta\big(|f(z)|^{p}\big)dV_{N}(z)<\infty$$
if and only if  $f\in b_{\alpha-2,p}(\mathbb{B}^{n})$.
\end{cor}

\begin{defn} For $m\in\{2,3,\ldots\}$, we denote by $\mathcal{HZ}_{m}(\mathbb{B}^{n})$ the class of all
functions $f\in\mathcal{C}^{m}(\mathbb{B}^{n})$ satisfying {\it
Heinz's} nonlinear differential inequality (cf. \cite{HZ})
$$|\Delta f(z)|\leq a(z)\|D_{f}(z)\|+b(z)|f(z)|+c(z),$$
where $a(z)$, $b(z)$ and $c(z)$ are real-valued nonnegative
continuous functions in $\mathbb{B}^{n}$.
\end{defn}

\begin{thm}\label{thm-1}
Let $M$ be a nonnegative constant and
$f\in\mathcal{HZ}_{3}(\mathbb{B}^{n})\cap\mathcal{D}_{\gamma,\alpha}(\mathbb{B}^{n})$
with ${\rm Re}(f\overline{\Delta f})\geq0$ and ${\rm
Re}\left\{\sum_{k=1}^{n}\big[\overline{f_{z_{k}}}(\Delta
f)_{z_{k}}+\overline{f_{\overline{z}_{k}}}(\Delta
f)_{\overline{z}_{k}}\big]\right\}\geq0,$ where $2\leq \alpha\leq
2n$, $\gamma>0$, $\sup_{z\in\mathbb{B}^{n}}a(z)<\infty$,
$\sup_{z\in\mathbb{B}^{n}}b(z)<\infty$ and $c(z)\leq
M\big(d(z)\big)^{-q}$. Then for $p\geq2$,
$$\int_{\mathbb{B}^{n}}\big(d(z)\big)^{pq}\Delta\big(|f(z)|^{p}\big)dV_{N}(z)<\infty,$$
where $q=(2n+\gamma)/\alpha-1.$
\end{thm}

The result given below is a consequence of Theorem \ref{thm-1}.

\begin{cor}\label{cor-3}
For $2\leq \alpha\leq 2n$ and $\gamma>0$, let
$f\in\mathcal{HZ}_{3}(\mathbb{B}^{n})\cap\mathcal{D}_{\gamma,\alpha}(\mathbb{B}^{n})$
satisfying {\rm (\ref{eq-CR1.1})}, where  $\tau$  and $\eta$ are
nonnegative constant functions. Then for $p\geq2$,
$$\int_{\mathbb{B}^{n}}\big(d(z)\big)^{pq}\Delta\big(|f(z)|^{p}\big)dV_{N}(z)<\infty,$$
where $q=(2n+\gamma)/\alpha-1.$
\end{cor}
\bpf By elementary calculations, we see that if $f$ is a solution to
(\ref{eq-CR1.1}), then $f$ satisfies Heinz's nonlinear differential
inequality. Hence Corollary \ref{cor-3} follows from
(\ref{eq-CR34}), (\ref{eq-CR37}) and Theorem \ref{thm-1}. \epf

By Corollaries \ref{cor-1}, \ref{cor-2} and \ref{cor-3}, we get

\begin{cor}\label{cor-4}
For $2\leq \alpha\leq 2n$ and $\gamma>0$, let
$q=(2n+\gamma)/\alpha-1$ and
$f\in\mathcal{HZ}_{3}(\mathbb{B}^{n})\cap\mathcal{D}_{\gamma,\alpha}(\mathbb{B}^{n})$
satisfying {\rm (\ref{eq-CR1.1})}, where  $\tau$  and $\eta$ are
nonnegative constant functions.\\
\item{{\rm(1)}} If
$p=\frac{1}{q}\geq2$, then  $f\in
\mathcal{H}_{g}^{p}(\mathbb{B}^{n})$;\\
\item{{\rm(2)}} If $p\geq2$ and
$pq\geq2$, then  $f\in b_{pq-2,p}(\mathbb{B}^{n})$.
\end{cor}

\begin{defn} For $p\geq2$, $t_{1}>0$, $t_{2}>0$ and $m\in\{2,3,\ldots\}$, we denote by
 $\mathcal{IHZ}_{m}^{t_{1},t_{2}}(\mathbb{B}^{n})$ the class of all
functions $f\in\mathcal{C}^{m}(\mathbb{B}^{n})$ satisfying the
inverse {\it Heinz's} nonlinear differential inequality
$$\Delta (|f(z)|^{p})\geq a_{1}(z)\|D_{f}(z)\|^{t_{1}}+b_{1}(z)|f(z)|^{t_{2}}+c_{1}(z),$$
where $a_{1}(z)$, $b_{1}(z)$ and $c_{1}(z)$ are real-valued
nonnegative continuous functions in $\mathbb{B}^{n}$.
\end{defn}

\begin{thm}\label{thm-4}
Let
$f\in\mathcal{IHZ}_{2}^{t_{1},t_{2}}(\mathbb{B}^{n})\cap\mathcal{H}_{g}^{p}(\mathbb{B}^{n})$,
where
$\inf_{z\in\mathbb{B}^{n}}a_{1}(z)+\inf_{z\in\mathbb{B}^{n}}b_{1}(z)>0$
and $\inf_{z\in\mathbb{B}^{n}}c_{1}(z)\geq0$.\\
\item{{\rm(1)}} If $\inf_{z\in\mathbb{B}^{n}}a_{1}(z)>0$, then
$f\in\mathcal{D}_{1,t_{1}}(\mathbb{B}^{n})$;\\
\item{{\rm(2)}} If $\inf_{z\in\mathbb{B}^{n}}b_{1}(z)>0$,
 then $ f\in
b_{1,t_{2}}(\mathbb{B}^{n})$.
\end{thm}

For $k\in\{1,\ldots,n\}$, let $\lambda_{k}\in\mathbb{R}$ be a
constant and $f\in\mathcal{C}^{1}(\mathbb{B}^{n})$ satisfying the
following nonlinear PDE, \be\label{eq-}\frac{\partial f}{\partial
\overline{z}_{k}}=\lambda_{k}|f|^{\alpha},\ee where
 $\alpha\geq 0.$ If, for each  $k\in\{1,\ldots,n\}$,
$\lambda_{k}=0$, then $f$ is holomorphic. Moreover, if $\alpha=0$,
then $f$ is pluriharmonic (cf. \cite{CR,Rudin}). It has attracted the
attention of many authors when $n=\lambda_{1}=1$ and
$\alpha\in(0,1)$ (cf. \cite{AIM,CY,IP}).

\begin{cor}\label{thm-4.0} For
$\sum_{k=1}^{n}\lambda_{k}^{2}\neq0$, $\alpha\geq 0$  and
$p>\max\{2, (\alpha-2)^{2}/4\}$, if $f\in
\mathcal{H}_{g}^{p}(\mathbb{B}^{n})\cap\mathcal{C}^{2}(\mathbb{B}^{n})$
satisfying {\rm (\ref{eq-})}, then $f\in
b_{1,\vartheta}(\mathbb{B}^{n})$, where $\vartheta=p+2\alpha-2.$
\end{cor}

The proofs of Theorem \ref{thm-3} will be presented in Section
\ref{csw-sec2}, and the proofs of Theorems \ref{thm-5.0},
\ref{thm-2}, \ref{thm-1}, \ref{thm-4}  and Corollary \ref{thm-4.0}
will be given in Section \ref{csw-sec3}.

\section{Hardy-type spaces and applications to pdes}\label{csw-sec2}
We start this section by recalling the following result. 

\begin{Thm}{\rm \cite{Pav}}\label{Green-thm}
Let $g$ be a function of class $\mathcal{C}^{2}(\mathbb{B}^{n})$.
Then, for $r\in(0,1)$,
$$\int_{\partial \mathbb{B}^{n}}g(r\zeta)\,d\sigma(\zeta)=
g(0)+\int_{ \mathbb{B}^{n}(r)}\Delta g(z)G_{2n}(z,r)\,dV_{N}(z),
$$
where $$ G_{2n}(z,r)=
\begin{cases}
\displaystyle\frac{|z|^{2(1-n)}-r^{2(1-n)}}{4n(n-1)},
& \mbox{ if } n\geq2,\\
\displaystyle\frac{1}{2}\log\frac{r}{|z|}, &\mbox{ if }\, n=1.
\end{cases}$$
\end{Thm}

\begin{lem}\label{lem-4} Let $p\geq2$ and $f\in\mathcal{C}^{2}(\mathbb{B}^{n})$ with ${\rm
Re}(f\overline{\Delta f})\geq0.$  Then  $M_{p}^{p}(r,f)$ is
increasing with respect to  $r\in(0,1)$.
\end{lem}
\bpf \noindent ${\rm \mathbf{Case~ 1.}}$ Let $p\in[4,\infty)$.

By elementary calculations, we get
$$\Delta\big(|f|^{p}\big)=p(p-2)|f|^{p-4}\sum_{k=1}^{n}|f_{z_{k}}\overline{f}+\overline{f_{\overline{z}_{k}}}f|^{2}+
2p|f|^{p-2}\|D_{f}\|^{2}+p|f|^{p-2}\mbox{Re}(f\overline{\Delta
f})\geq0,$$ which implies that, for $p\in[4,\infty)$,
$M_{p}^{p}(r,f)$ is increasing in  $(0,1)$.

\noindent ${\rm \mathbf{Case~ 2.}}$ Let $p\in[2,4)$.

 For
$m\in\{1,2,\ldots\}$, let
$T_{m}^{p}=\left(|f|^{2}+\frac{1}{m}\right)^{\frac{p}{2}}$. By
computations, we have
\begin{eqnarray*}
\Delta\big(T_{m}^{p}\big)&=&4\sum_{k=1}^{n}\frac{\partial^{2}}{\partial
z_{k}\partial\overline{z}_{k}}
\big(T_{m}^{p}\big)=4\sum_{k=1}^{n}\big(T_{m}^{p}\big)_{z_{k}\overline{z}_{k}}\\
&=&p(p-2)\left(|f|^{2}+\frac{1}{m}\right)^{\frac{p}{2}-2}\sum_{k=1}^{n}|f_{z_{k}}\overline{f}+\overline{f_{\overline{z}_{k}}}f|^{2}\\
&&+2p\left(|f|^{2}+\frac{1}{m}\right)^{\frac{p}{2}-1}\|D_{f}\|^{2}+p\left(|f|^{2}+\frac{1}{m}\right)^{\frac{p}{2}-1}\mbox{Re}(f\overline{\Delta
f}).
\end{eqnarray*}
Let $Q_{m}=\Delta\big(T_{m}^{p}\big)$. It is not difficult to show
that, for $r\in(0,1)$, $Q_{m}$ is integrable in $\mathbb{B}^{n}(r)$
and $0<Q_{m}\leq \Lambda_{f},$ where
\begin{eqnarray*}
\Lambda_{f}&=&p(p-2)|f|^{p-2}
\sum_{k=1}^{n}\big(|f_{z_{k}}|+|f_{\overline{z}_{k}}|\big)^{2}
+2p\left(|f|^{2}+1\right)^{\frac{p}{2}-1}\|D_{f}\|^{2}\\
&&+p\left(|f|^{2}+1\right)^{\frac{p}{2}-1}\mbox{Re}(f\overline{\Delta
f})
\end{eqnarray*}
and $\Lambda_{f}$ is integrable in $\mathbb{B}^{n}(r)$.

By using Theorem \Ref{Green-thm} and Lebesgue's Dominated Convergence
theorem, 
we  get

\begin{eqnarray*}
\lim_{m\rightarrow\infty}r^{2n-1}\frac{d}{dr}M_{p}^{p}(r,T_{m})&=&\frac{1}{2n}
\lim_{m\rightarrow\infty}\int_{\mathbb{B}^{n}(r)}Q_{m}dV_{N}\\
&=&\frac{1}{2n}\int_{\mathbb{B}^{n}(r)}\lim_{m\rightarrow\infty}Q_{m}dV_{N}\\
&=&\frac{1}{2n}\int_{\mathbb{B}^{n}(r)}\big[p(p-2)|f|^{p-4}\sum_{k=1}^{n}|f\overline{f_{z_{k}}}+\overline{f}f_{\overline{z}_{k}}|^{2}\\
&&+2p|f|^{p-2}\|D_{f}\|^{2}+p|f|^{p-2}\mbox{Re}\big(f\overline{\Delta f}\big)\big]dV_{N}\\
 &=&r^{2n-1}\frac{d}{dr}M_{p}^{p}(r,f)\\
 &\geq&0,
\end{eqnarray*}
which implies that $M_{p}^{p}(r,f)$ is increasing in $r$ on $(0,1)$  for
$p\in[2,4)$.

 \epf

By using Theorem \Ref{Green-thm} and a similar argument as in the
proof of Lemma \ref{lem-4}, we obtain the following result.

\begin{lem}\label{lem-5} Let $p\geq2$ and $f\in\mathcal{C}^{2}(\mathbb{B}^{n})$ with ${\rm
Re}(f\overline{\Delta f})\geq0.$  Then,  for $r\in(0,1)$,
$$M_{p}^{p}(r,f)=|f(0)|^{p}+ \int_{\mathbb{B}^{n}(r)}\Delta
(|f(z)|^{p})G_{2n}(z,r)dV_{N}(z),$$  where $G_{2n}$ is the function
defined in Theorem {\rm \Ref{Green-thm}}.
\end{lem}

\subsection*{Proof of Theorem \ref{thm-3}}
\noindent ${\rm \mathbf{Case~ 1.}}$ Let $n\geq2$.

We first prove the necessity. For a fixed positive constant
$r_{0}\in(0,1)$, let $r\in(r_{0},1)$. Then,  by Lemma \ref{lem-5},
we have

\beq\label{eq-CR19} \nonumber
M_{p}^{p}(r,f)&=&|f(0)|^{p}+\int_{\mathbb{B}^{n}(r)}\Delta\big(|f(z)|^{p}\big)G_{2n}(z,r)dV_{N}(z)\\
\nonumber&=&|f(0)|^{p}+\frac{1}{4n(n-1)}\int_{\mathbb{B}^{n}(r)\setminus\mathbb{B}^{n}(r_{0})}
\big(|z|^{2(1-n)}-r^{2(1-n)}\big)\Delta\big(|f(z)|^{p}\big)dV_{N}(z)\\
 &&+
\frac{1}{4n(n-1)}\int_{\mathbb{B}^{n}(r_{0})}\big(|z|^{2(1-n)}-r^{2(1-n)}\big)\Delta\big(|f(z)|^{p}\big)dV_{N}(z).\eeq
Since $\Delta\big(|f|^{p}\big)\geq0$,

\beq\label{eq-CR20}
\infty&>&2n\int_{\partial\mathbb{B}^{n}}\int_{0}^{r_{0}}\big(\rho-
\rho^{2n-1}\big)\Delta\big(|f(\rho\zeta)|^{p}\big)d\rho
d\sigma(\zeta)\\  \nonumber
&\geq&2n\int_{\partial\mathbb{B}^{n}}\int_{0}^{r_{0}}\big(\rho-r^{2(1-n)}
\rho^{2n-1}\big)\Delta\big(|f(\rho\zeta)|^{p}\big)d\rho
d\sigma(\zeta)\\ \nonumber
&=&\int_{\mathbb{B}^{n}(r_{0})}\big(|z|^{2(1-n)}-r^{2(1-n)}\big)\Delta\big(|f(z)|^{p}\big)dV_{N}(z)
\eeq and

\beq\label{eq-CR21}
\nonumber&&\int_{\mathbb{B}^{n}(r)\setminus\mathbb{B}^{n}(r_{0})}\big(|z|^{2(1-n)}-r^{2(1-n)}\big)\Delta\big(|f(z)|^{p}\big)dV_{N}(z)\\
\nonumber
&=&\int_{\mathbb{B}^{n}(r)\setminus\mathbb{B}^{n}(r_{0})}\frac{(r-|z|)\big(\sum_{k=0}^{2n-3}r^{2n-3-k}|z|^{k}\big)}{|z|^{2n-2}r^{2n-2}}
\Delta\big(|f(z)|^{p}\big)dV_{N}(z)\\
 \nonumber &\leq&
\frac{(2n-2)}{r_{0}^{4n-4}}\int_{\mathbb{B}^{n}(r)\setminus\mathbb{B}^{n}(r_{0})}(r-|z|)\Delta\big(|f(z)|^{p}\big)dV_{N}(z)\\
\nonumber
&\leq&\frac{(2n-2)}{r_{0}^{4n-4}}\int_{\mathbb{B}^{n}}d(z)\Delta\big(|f(z)|^{p}\big)dV_{N}(z)\\
\nonumber &<&\infty.
 \eeq
By  (\ref{eq-CR19}) and Lemma \ref{lem-4}, we see that the limit
$$\lim_{r\rightarrow1-}M_{p}(r,f)$$ exists. Hence
$f\in \mathcal{H}_{g}^{p}(\mathbb{B}^{n})$.

Next we prove that sufficiency. Applying (\ref{eq-CR19}),
(\ref{eq-CR20}) and $f\in \mathcal{H}_{g}^{p}(\mathbb{B}^{n})$, we
observe that
 \beq\label{eq-CR22}
\nonumber\infty&>&\int_{\mathbb{B}^{n}(r)\setminus\mathbb{B}^{n}(r_{0})}\big(|z|^{2(1-n)}-r^{2(1-n)}\big)\Delta\big(|f(z)|^{p}\big)dV_{N}(z)\\
\nonumber
&=&\int_{\mathbb{B}^{n}(r)\setminus\mathbb{B}^{n}(r_{0})}\frac{(r-|z|)\big(\sum_{k=0}^{2n-3}r^{2n-3-k}|z|^{k}\big)}{|z|^{2n-2}r^{2n-2}}
\Delta\big(|f(z)|^{p}\big)dV_{N}(z)\\ \nonumber &\geq&I(r), \eeq
which, together with the monotonicity of $I(r)$ on $r\in[r_{0},1)$,
yields that
$$\lim_{r\rightarrow1-}\int_{\mathbb{B}^{n}(r)\setminus\mathbb{B}^{n}(r_{0})}(r-|z|)\Delta\big(|f(z)|^{p}\big)dV_{N}(z)$$
exists, where
$$I(r)=(2n-2)r_{0}^{2n-3}\int_{\mathbb{B}^{n}(r)\setminus\mathbb{B}^{n}(r_{0})}(r-|z|)\Delta\big(|f(z)|^{p}\big)dV_{N}(z).$$
Therefore,
$$\int_{\mathbb{B}^{n}}d(z)\Delta\big(|f(z)|^{p}\big)dV_{N}(z)<\infty.$$

\noindent ${\rm \mathbf{Case~ 2.}}$ Let $n=1$.

In this case, we also first prove  necessity. 
Fix $r\in(0,1)$. Since
$$\lim_{|z|\rightarrow r}\frac{\log r-\log
|z|}{r-|z|}=\frac{1}{r},$$ we see that there exists $r_{0}\in(0,r)$ such
that \beq\label{eq-CR25}\frac{1}{2r}\leq \frac{\log r-\log
|z|}{r-|z|}\leq\frac{3}{2r}\eeq for $r_{0}\leq|z|<r.$ It is not
difficult to see that, for $|z|\leq r<1,$ \be\label{eq-0.1}
\frac{r-|z|}{r}\leq1-|z|.\ee Because
 $$\lim_{\rho\rightarrow0+}\rho\log\frac{1}{\rho}=0,$$ it follows that

\beq\label{eq-CR23}
\int_{\mathbb{D}(r_{0})}\Delta\big(|f(z)|^{p}\big)\log\frac{r}{|z|}dA(z)&=&
\int_{\mathbb{D}(r_{0})}\Delta\big(|f(z)|^{p}\big)\log\frac{1}{|z|}dA(z)\\
\nonumber
&=&\frac{1}{\pi}\int_{0}^{2\pi}\int_{0}^{r_{0}}\Delta\big(|f(\rho
e^{i\theta})|^{p}\big)\rho\log\frac{1}{\rho}d\theta\\
\nonumber&<&\infty,\eeq where
 $dA$ denotes the
normalized area measure in $\mathbb{D}$.

By (\ref{eq-CR25}), (\ref{eq-0.1}),  (\ref{eq-CR23}), Lemmas
\ref{lem-4} and \ref{lem-5}, we see that

\begin{eqnarray*}
M_{p}^{p}(r,f)& =&|f(0)|^{p}+
\frac{1}{2}\int_{\mathbb{D}(r)}\Delta
\big(|f(z)|^{p}\big)\log\frac{r}{|z|}\,dA(z)\\
&=&|f(0)|^{p}+ \frac{1}{2}\int_{\mathbb{D}(r_{0})}\Delta
\big(|f(z)|^{p}\big)\log\frac{r}{|z|}\,dA(z)\\
&&+\frac{1}{2}\int_{\mathbb{D}(r)\setminus\mathbb{D}(r_{0})}\Delta
\big(|f(z)|^{p}\big)\log\frac{r}{|z|}\,dA(z)\\
&\leq& |f(0)|^{p}+ \frac{1}{2}\int_{\mathbb{D}(r_{0})}\Delta
\big(|f(z)|^{p}\big)\log\frac{r}{|z|}\,dA(z)\\
&&+\frac{3}{4}\int_{\mathbb{D}(r)\setminus\mathbb{D}(r_{0})}\Delta
\big(|f(z)|^{p}\big)\frac{(r-|z|)}{r}\,dA(z)\\
 &\leq& |f(0)|^{p}+
\frac{1}{2}\int_{\mathbb{D}(r_{0})}\Delta
\big(|f(z)|^{p}\big)\log\frac{r}{|z|}\,dA(z)\\
&&+\frac{3}{4}\int_{\mathbb{D}\setminus\mathbb{D}(r_{0})}\Delta
\big(|f(z)|^{p}\big)d(z)\,dA(z)\\
 &<&\infty,
\end{eqnarray*} which implies that the limit
$$\lim_{r\rightarrow1-}M_{p}(r,f)$$  exists. Hence
$f\in\mathcal{H}_{g}^{p}(\mathbb{B}^{n}).$

Now we prove that sufficiency. By (\ref{eq-CR25}), we have
\beq\label{eq-CR24}
\nonumber\int_{\mathbb{D}(r)\setminus\mathbb{D}(r_{0})}\Delta
\big(|f(z)|^{p}\big)\log\frac{r}{|z|}\,dA(z)&\geq&\frac{1}{2r}\int_{\mathbb{D}(r)\setminus\mathbb{D}(r_{0})}\Delta
\big(|f(z)|^{p}\big)(r-|z|)\,dA(z)\\
&\geq&\frac{I^{\ast}(r)}{2}, \eeq where
$$I^{\ast}(r)=\int_{\mathbb{D}(r)\setminus\mathbb{D}(r_{0})}\Delta
\big(|f(z)|^{p}\big)(r-|z|)\,dA(z).$$

 By (\ref{eq-CR24}), Lemmas \ref{lem-4} and
\ref{lem-5}, we have
\begin{eqnarray*}
M_{p}^{p}(r,f)& =&|f(0)|^{p}+ \frac{1}{2}\int_{\mathbb{D}(r)}\Delta
\big(|f(z)|^{p}\big)\log\frac{r}{|z|}\,dA(z)\\
&=&|f(0)|^{p}+ \frac{1}{2}\int_{\mathbb{D}(r_{0})}\Delta
\big(|f(z)|^{p}\big)\log\frac{r}{|z|}\,dA(z)\\
&&+\frac{1}{2}\int_{\mathbb{D}(r)\setminus\mathbb{D}(r_{0})}\Delta
\big(|f(z)|^{p}\big)\log\frac{r}{|z|}\,dA(z)\\
&\geq&|f(0)|^{p}+ \frac{1}{2}\int_{\mathbb{D}(r_{0})}\Delta
\big(|f(z)|^{p}\big)\log\frac{r}{|z|}\,dA(z)+\frac{1}{4}I^{\ast}(r),
\end{eqnarray*}
which yields that $I^{\ast}(r)<\infty.$ Since $I^{\ast}(r)$ is
increasing on $r$, we see that $$\lim_{r\rightarrow1-}I^{\ast}(r)$$
 exists. Then
$$\int_{\mathbb{D}}d(z)\Delta\big(|f(z)|^{p}\big)dA(z)<\infty$$
concluding the proof of the theorem. \qed

\begin{lem}\label{lem-1}
Let $f\in\mathcal{C}^{3}(\mathbb{B}^{n})$ and ${\rm
Re}\left\{\sum_{k=1}^{n}\big[\overline{f_{z_{k}}}(\Delta
f)_{z_{k}}+\overline{f_{\overline{z}_{k}}}(\Delta
f)_{\overline{z}_{k}}\big]\right\}\geq0.$ Then,  for $\alpha\geq2$,
$\|D_{f}\|^{\alpha}$ is subharmonic in $\mathbb{B}^{n}$.
\end{lem}
\bpf First we consider  the case $\alpha\in[4,\infty)$. Since
\begin{eqnarray*}
\Delta\big(\|D_{f}\|^{\alpha}\big)&=&\alpha(\alpha-2)\|D_{f}\|^{\alpha-4}\left|\sum_{j=1}^{n}\sum_{k=1}^{n}
\big(f_{z_{k}z_{j}}\overline{f_{z_{k}}}+\overline{f_{z_{k}\overline{z}_{j}}}f_{z_{k}}+
f_{\overline{z}_{k}z_{j}}\overline{f_{\overline{z}_{k}}}+\overline{f_{\overline{z}_{k}\overline{z}_{j}}}f_{\overline{z}_{k}}\big)\right|^{2}\\
&&+2\alpha\|D_{f}\|^{\alpha-2}\sum_{j=1}^{n}\sum_{k=1}^{n}
\big(|f_{z_{j}z_{k}}|^{2}+|f_{z_{j}\overline{z}_{k}}|^{2}+|f_{\overline{z}_{j}z_{k}}|^{2}+|f_{\overline{z}_{j}\overline{z}_{k}}|^{2}\big)\\
&&+\alpha\|D_{f}\|^{\alpha-2}{\rm
Re}\left\{\sum_{k=1}^{n}\Big[\overline{f_{z_{k}}}(\Delta
f)_{z_{k}}+\overline{f_{\overline{z}_{k}}}(\Delta
f)_{\overline{z}_{k}}\Big]\right\} \geq0,
\end{eqnarray*} we see that, for $\alpha\in[4,\infty)$, $\|D_{f}\|^{\alpha}$ is subharmonic in $\mathbb{B}^{n}$.

 Next we deal with the case $\alpha\in[2,4).$ In this case, for
$m\in\{1,2,\ldots\}$, we let
$F_{m}^{\alpha}=(\|D_{f}\|^{2}+\frac{1}{m})^{\frac{\alpha}{2}}.$
Then, by elementary computations, we have

\begin{eqnarray*}
\Delta(F_{m}^{\alpha})&=&4\sum_{j=1}^{n}(F_{m}^{\alpha})_{z_{k}\overline{z}_{k}}\\
&=&4\sum_{j=1}^{n}\frac{\partial^{2}}{\partial z_{j}\partial
\overline{z}_{j}}\Bigg\{\left[\frac{1}{m}+\sum_{k=1}^{n}\big(f_{z_{k}}\overline{f_{z_{k}}}+f_{\overline{z}_{k}}\overline{f_{\overline{z}_{k}}}\big)
\right]^{\frac{\alpha}{2}}\Bigg\}\\
&=&\alpha(\alpha-2)\left(\|D_{f}\|^{2}+\frac{1}{m}\right)^{\frac{\alpha}{2}-2}\left[\sum_{j=1}^{n}\frac{\partial}{\partial
z_{j}}(\|D_{f}\|^{2})\right]\left[\sum_{j=1}^{n}\frac{\partial}{\partial
\overline{z}_{j}}(\|D_{f}\|^{2})\right]\\
&&+2\alpha\left(\|D_{f}\|^{2}+\frac{1}{m}\right)^{\frac{\alpha}{2}-1}\sum_{j=1}^{n}\sum_{k=1}^{n}
\big(|f_{z_{j}z_{k}}|^{2}+|f_{z_{j}\overline{z}_{k}}|^{2}+|f_{\overline{z}_{j}z_{k}}|^{2}+|f_{\overline{z}_{j}\overline{z}_{k}}|^{2}\big)\\
&&+\alpha\left(\|D_{f}\|^{2}+\frac{1}{m}\right)^{\frac{\alpha}{2}-1}{\rm
Re}\left\{\sum_{k=1}^{n}\big[\overline{f_{z_{k}}}(\Delta
f)_{z_{k}}+\overline{f_{\overline{z}_{k}}}(\Delta
f)_{\overline{z}_{k}}\big]\right\}\\
&=&\alpha(\alpha-2)\left(\|D_{f}\|^{2}+\frac{1}{m}\right)^{\frac{\alpha}{2}-2}\\
&&\times\left|\sum_{j=1}^{n}\sum_{k=1}^{n}
\big(f_{z_{k}z_{j}}\overline{f_{z_{k}}}+\overline{f_{z_{k}\overline{z}_{j}}}f_{z_{k}}+
f_{\overline{z}_{k}z_{j}}\overline{f_{\overline{z}_{k}}}+\overline{f_{\overline{z}_{k}\overline{z}_{j}}}f_{\overline{z}_{k}}\big)\right|^{2}\\
&&+2\alpha\left(\|D_{f}\|^{2}+\frac{1}{m}\right)^{\frac{\alpha}{2}-1}\sum_{j=1}^{n}\sum_{k=1}^{n}
\big(|f_{z_{j}z_{k}}|^{2}+|f_{z_{j}\overline{z}_{k}}|^{2}+|f_{\overline{z}_{j}z_{k}}|^{2}+|f_{\overline{z}_{j}\overline{z}_{k}}|^{2}\big)\\
&&+\alpha\left(\|D_{f}\|^{2}+\frac{1}{m}\right)^{\frac{\alpha}{2}-1}{\rm
Re}\left\{\sum_{k=1}^{n}\Big[\overline{f_{z_{k}}}(\Delta
f)_{z_{k}}+\overline{f_{\overline{z}_{k}}}(\Delta
f)_{\overline{z}_{k}}\Big]\right\}.
\end{eqnarray*}

By the Cauchy-Schwarz inequality, we have

\beq\label{eq-CR33} &&\bigg|\sum_{j=1}^{n}\sum_{k=1}^{n}
\Big(f_{z_{k}z_{j}}\overline{f_{z_{k}}}+\overline{f_{z_{k}\overline{z}_{j}}}f_{z_{k}}+
f_{\overline{z}_{k}z_{j}}\overline{f_{\overline{z}_{k}}}+
\overline{f_{\overline{z}_{k}\overline{z}_{j}}}f_{\overline{z}_{k}}\Big)\bigg|^{2}\\
\nonumber &\leq&\bigg[\sum_{j=1}^{n}\sum_{k=1}^{n}
\Big(|f_{z_{k}z_{j}}\overline{f_{z_{k}}}|+|\overline{f_{z_{k}\overline{z}_{j}}}f_{z_{k}}|+
|f_{\overline{z}_{k}z_{j}}\overline{f_{\overline{z}_{k}}}|+
|\overline{f_{\overline{z}_{k}\overline{z}_{j}}}f_{\overline{z}_{k}}|\Big)\bigg]^{2}\\
\nonumber
&\leq&\Bigg\{\sum_{j=1}^{n}\sum_{k=1}^{n}\left[\Big(2|f_{z_{k}}|^{2}+2|f_{\overline{z}_{k}}|^{2}\Big)^{\frac{1}{2}}
\Big(|f_{z_{k}z_{j}}|^{2}+|f_{z_{k}\overline{z}_{j}}|^{2}+
|f_{\overline{z}_{k}z_{j}}|^{2}+
|f_{\overline{z}_{k}\overline{z}_{j}}|^{2}\Big)^{\frac{1}{2}}\right]\Bigg\}^{2}\\
\nonumber
&\leq&2\sum_{j=1}^{n}\sum_{k=1}^{n}\left[\Big(2|f_{z_{k}}|^{2}+2|f_{\overline{z}_{k}}|^{2}\Big)
\Big(|f_{z_{k}z_{j}}|^{2}+|f_{z_{k}\overline{z}_{j}}|^{2}+
|f_{\overline{z}_{k}z_{j}}|^{2}+
|f_{\overline{z}_{k}\overline{z}_{j}}|^{2}\Big)\right]\\ \nonumber
&\leq&4\|D_{f}\|^{2}\sum_{j=1}^{n}\sum_{k=1}^{n}
\Big(|f_{z_{k}z_{j}}|^{2}+|f_{z_{k}\overline{z}_{j}}|^{2}+
|f_{\overline{z}_{k}z_{j}}|^{2}+
|f_{\overline{z}_{k}\overline{z}_{j}}|^{2}\Big). \eeq

  Hence, by (\ref{eq-CR33}) and Lebesgue's Dominated
Convergence theorem, we obtain
\begin{eqnarray*}
\lim_{m\rightarrow\infty}\Delta(F_{m}^{\alpha})&=&\alpha(\alpha-2)\|D_{f}\|^{\alpha-4}\left|\sum_{j=1}^{n}\sum_{k=1}^{n}
\big(f_{z_{k}z_{j}}\overline{f_{z_{k}}}+\overline{f_{z_{k}\overline{z}_{j}}}f_{z_{k}}+
f_{\overline{z}_{k}z_{j}}\overline{f_{\overline{z}_{k}}}+\overline{f_{\overline{z}_{k}\overline{z}_{j}}}f_{\overline{z}_{k}}\big)\right|^{2}\\
&&+2\alpha\|D_{f}\|^{\alpha-2}\sum_{j=1}^{n}\sum_{k=1}^{n}
\big(|f_{z_{j}z_{k}}|^{2}+|f_{z_{j}\overline{z}_{k}}|^{2}+|f_{\overline{z}_{j}z_{k}}|^{2}+|f_{\overline{z}_{j}\overline{z}_{k}}|^{2}\big)\\
&&+\alpha\|D_{f}\|^{\alpha-2}{\rm
Re}\left\{\sum_{k=1}^{n}\Big[\overline{f_{z_{k}}}(\Delta
f)_{z_{k}}+\overline{f_{\overline{z}_{k}}}(\Delta
f)_{\overline{z}_{k}}\Big]\right\} \geq0.
\end{eqnarray*} Then, for $\alpha\in[2,4)$,
$\|D_{f}\|^{\alpha}$ is subharmonic in $\mathbb{B}^{n}$.
 \epf

\subsection*{Proof of Theorem
\ref{thm-5.0}} It is not difficult to see that if  $\tau$ and
$\eta$ are
 constant functions, then each solution $f$ to
(\ref{eq-CR1.1}) belongs to $C^{\infty}(\mathbb{B}^{n})$, i.e., they
are infinitely differentiable in $\mathbb{B}^{n}$.

By elementary calculations, we get

\beq\label{eq-CR34} &&\sum_{k=1}^{n}{\rm
Re}\Big[\overline{f_{z_{k}}(z)}(\Delta
f(z))_{z_{k}}+\overline{f_{\overline{z}_{k}}(z)}(\Delta
f(z))_{\overline{z}_{k}}\Big]\\ \nonumber
&=&\sum_{k=1}^{n}\left[\tau\big(|f_{z_{k}}|^{2}+
|f_{\overline{z}_{k}}|^{2}\big)+\frac{\eta}{2}\big|f_{z_{k}}+\overline{f_{\overline{z}_{k}}}\big|^{2}\right]\\
\nonumber &\leq&(\eta+\tau)\|D_{f}\|^{2}, \eeq and
\beq\label{eq-CR37} {\rm Re}\big(f\overline{\Delta
f}\big)=\tau|f|^{2}+\eta\big(\mbox{Re}(f)\big)^{2}\leq
(\eta+\tau)|f|^{2}. \eeq

By using H\"older's inequality, for $\rho\in(0,1),$ we see that

 \beq\label{eq-CR35}
\int_{\partial\mathbb{B}^{n}}\|D_{f}(\rho\zeta)\|^{p-2}
\big(D^{\ast}_{f}(\rho\zeta)\big)^{2}d\sigma(\zeta)\leq
M_{p}^{2}(\rho,D^{\ast}_{f})M_{p}^{p-2}(\rho,D_{f})
 \eeq and

\beq\label{eq-CR36}
\int_{\partial\mathbb{B}^{n}}\|f(\rho\zeta)\|^{p-2}
\|D_{f}(\rho\zeta)\|^{2}d\sigma(\zeta)\leq
M_{p}^{2}(\rho,D_{f})M_{p}^{p-2}(\rho,f).
 \eeq

For $t\in[0,1]$, $r\in(0,1)$ and $\rho\in(0,r]$, we obtain

\be\label{eq-CR-35} \frac{t(1-t^{2n-2})}{2(n-1)}\leq1-t \ee and
\be\label{eq-CR-36} \rho\log\frac{r}{\rho}\leq r-\rho, \ee where
$n\geq2.$

 \noindent ${\rm
\mathbf{Case~ 1.}}$ Let $n\geq2$.

\noindent ${\rm \mathbf{Step~ 1.}}$ By (\ref{eq-CR33}),
(\ref{eq-CR34}), (\ref{eq-CR35}), Lemma \ref{lem-1}, Theorem
\Ref{Green-thm} and Lebesgue's Dominated Convergence theorem, we see
that

\begin{eqnarray*}
M_{p}^{p}(r,D_{f})&=&\|D_{f}(0)\|^{p}+\int_{\mathbb{B}^{n}(r)}\Delta\big(\|D_{f}(z)\|^{p}\big)G_{2n}(z,r)dV_{N}(z)\\
&=&\|D_{f}(0)\|^{p}\\
&&+\int_{\mathbb{B}^{n}(r)}\Bigg\{p\|D_{f}(z)\|^{p-2}\sum_{k=1}^{n}{\rm
Re}\Big[\overline{f_{z_{k}}(z)}(\Delta
f(z))_{z_{k}}+\overline{f_{\overline{z}_{k}}(z)}(\Delta
f(z))_{\overline{z}_{k}}\Big]\\
&&+p(p-2)\|D_{f}(z)\|^{p-4}\bigg|\sum_{j=1}^{n}\sum_{k=1}^{n}
\Big(f_{z_{k}z_{j}}(z)\overline{f_{z_{k}}(z)}\\
&&+\overline{f_{z_{k}\overline{z}_{j}}(z)}f_{z_{k}}(z)+
f_{\overline{z}_{k}z_{j}}(z)\overline{f_{\overline{z}_{k}}(z)}+
\overline{f_{\overline{z}_{k}\overline{z}_{j}}(z)}f_{\overline{z}_{k}}(z)\Big)\bigg|^{2}\\
&&+2p\|D_{f}(z)\|^{p-2}\big(D^{\ast}_{f}(z)\big)^{2}\Bigg\}G_{2n}(z,r)dV_{N}(z)\\
&\leq&\|D_{f}(0)\|^{p}+p\int_{\mathbb{B}^{n}(r)}\Big[(\eta+\tau)\|D_{f}(z)\|^{p}\\
&&+2(2p-3)\|D_{f}(z)\|^{p-2}\big(D^{\ast}_{f}(z)\big)^{2}
\Big]G_{2n}(z,r)dV_{N}(z)\\
&=&\|D_{f}(0)\|^{p}+\frac{p(\eta+\tau)}{2(n-1)}\int_{0}^{r}\big(\rho-\rho^{2n-1}r^{2(1-n)}\big)M_{p}^{p}(\rho,
D_{f})d\rho\\
&&+\frac{p(2p-3)}{(n-1)}\int_{0}^{r}\big(\rho-\rho^{2n-1}r^{2(1-n)}\big)\int_{\partial\mathbb{B}^{n}}\|D_{f}(\rho\zeta)\|^{p-2}
\big(D^{\ast}_{f}(\rho\zeta)\big)^{2}d\sigma(\zeta)d\rho\\
&\leq&\|D_{f}(0)\|^{p}+\frac{p(\eta+\tau)}{2(n-1)}\int_{0}^{r}\big(\rho-\rho^{2n-1}r^{2(1-n)}\big)M_{p}^{p}(\rho,
D_{f})d\rho\\
&&+\frac{p(2p-3)}{(n-1)}\int_{0}^{r}\big(\rho-\rho^{2n-1}r^{2(1-n)}\big)M_{p}^{2}(\rho,D^{\ast}_{f})M_{p}^{p-2}(\rho,D_{f})d\rho.
\end{eqnarray*}
The above together with (\ref{eq-CR-35}) and subharmonicity of
$\|D_{f}\|^{p}$,
 shows that

\begin{eqnarray*}
&&\left[1-\frac{p(\eta+\tau)}{2(n-1)}\int_{0}^{r}(\rho-\rho^{2n-1}r^{2(1-n)})d\rho\right]M_{p}^{2}(r,D_{f})\\
&=&\left[1-\frac{pr^{2}(\eta+\tau)}{4n}\right]M_{p}^{2}(r,D_{f})\\
&\leq&\|D_{f}(0)\|^{2}+\frac{p(2p-3)}{(n-1)}\int_{0}^{r}\big(\rho-\rho^{2n-1}r^{2(1-n)}\big)M_{p}^{2}(\rho,D^{\ast}_{f})d\rho\\
&=&\|D_{f}(0)\|^{2}+2p(2p-3)r^{2}\int_{0}^{1}\frac{t(1-t^{2n-2})}{2(n-1)}M_{p}^{2}\big(rt,D^{\ast}_{f}\big)dt\\
\end{eqnarray*}
\begin{eqnarray*}
&\leq&\|D_{f}(0)\|^{2}+2p(2p-3)r^{2}(M^{\ast})^{2}\int_{0}^{1}\left[\omega\Big(\frac{1}{1-rt}\Big)\right]^{2}(1-t)dt\\
&\leq&\|D_{f}(0)\|^{2}+2p(2p-3)r^{2}(M^{\ast})^{2}\int_{0}^{1}\left[\omega\Big(\frac{1}{1-rt}\Big)\right]^{2}(1-rt)dt\\
&\leq&\|D_{f}(0)\|^{2}+2p(2p-3)r^{2}(M^{\ast})^{2}\omega(1)\int_{0}^{1}\omega\Big(\frac{1}{1-rt}\Big)dt.
\end{eqnarray*}
Then \be\label{eq-CR38}M_{p}^{2}(r,D_{f})\leq M^{\ast}_{2}
\left[\|D_{f}(0)\|^{2}+M^{\ast}_{1}\int_{0}^{1}\omega\Big(\frac{1}{1-rt}\Big)dt\right],
\ee where $M^{\ast}_{1}=2p(2p-3)(M^{\ast})^{2}\omega(1)$ and
$M^{\ast}_{2}=1/\left[1-p(\eta+\tau)/(4n)\right].$

\noindent ${\rm \mathbf{Step~ 2.}}$ By (\ref{eq-CR37}),
(\ref{eq-CR36}), Lemmas \ref{lem-4} and \ref{lem-5}, we obtain

\begin{eqnarray*}
M_{p}^{p}(r, f)&=&|f(0)|^{p}+\int_{\mathbb{B}^{n}(r)}\Delta(|f(z)|^{p})G_{2n}(z,r)\,dV_{N}(z)\\
&\leq&|f(0)|^{p}+\int_{0}^{r}\int_{\partial\mathbb{B}^{n}}4np(p-1)\rho^{2n-1}|f(\rho\zeta)|^{p-2}
|D_{f}(\rho\zeta)|^{2}G_{2n}(\rho \zeta,r)d\sigma(\zeta)d\rho\\
&&+p(\eta+\tau)\int_{\mathbb{B}^{n}(r)}|f(z)|^{p}G_{2n}(z,r)\,dV_{N}(z)\\
&\leq&|f(0)|^{p}+4p(p-1)\int_{0}^{r}n\rho^{2n-1}G_{2n}(\rho\zeta,r)M_{p}^{2}(\rho,D_{f}) M_{p}^{p-2}(\rho, f)\,d\rho\\
&&+\frac{p(\eta+\tau) r^{2}}{4n}M_{p}^{p}(r,f).
\end{eqnarray*}
By the above estimates, (\ref{eq-CR-35}), (\ref{eq-CR38}) and the
monotonicity of $M_p(r,f)$ on $r$,
\begin{eqnarray*}
 \frac{M_{p}^{2}(r,
 f)}{M^{\ast}_{2}}
 &\leq&\left[1-\frac{p(\eta+\tau) r^{2}}{4n}\right]M_{p}^{2}(r,
 f)\\
 &\leq&|f(0)|^{2}+4p(p-1)\int_{0}^{r}n\rho^{2n-1}G_{2n}(\rho\zeta,r)M_{p}^{2}(\rho,D_{f})\,d\rho\\
&=&|f(0)|^{2}+2p(p-1)\int_{0}^{1}r^{2} M_{p}^{2}(r\rho,D_{f})\cdot\frac{\rho(1-\rho^{2n-2})}{2(n-1)}\,d\rho\\
&\leq&|f(0)|^{2}+2p(p-1)\int_{0}^{1}M_{p}^{2}(r\rho,D_{f})(1-\rho)\,d\rho\\
&\leq&|f(0)|^{2}+2p(p-1)M^{\ast}_{2}\|D_{f}(0)\|^{2}\\
&&+2p(p-1)M^{\ast}_{1}M^{\ast}_{2}\int_{0}^{1}\left[\int_{0}^{1}\omega\Big(\frac{1}{1-r\rho
t}\Big)(1-\rho)dt\right]d\rho\\
\end{eqnarray*}
\begin{eqnarray*}
&\leq&|f(0)|^{2}+2p(p-1)M^{\ast}_{2}\|D_{f}(0)\|^{2}\\
&&+2p(p-1)M^{\ast}_{1}M^{\ast}_{2}\int_{0}^{1}\left[\int_{0}^{1}\omega\Big(\frac{1}{1-r\rho
t}\Big)(1-rt\rho)dt\right]d\rho\\
&\leq&|f(0)|^{2}+2p(p-1)M^{\ast}_{2}\|D_{f}(0)\|^{2}+2p(p-1)M^{\ast}_{1}M^{\ast}_{2}\omega(1)\\
&<&\infty.
\end{eqnarray*}
Hence $f\in \mathcal{H}_{g}^{p}(\mathbb{B}^{n})$.

\noindent ${\rm \mathbf{Case~ 2.}}$ Let $n=1$.

\noindent ${\rm \mathbf{Step~ 3.}}$ By (\ref{eq-CR33}),
(\ref{eq-CR34}), (\ref{eq-CR35}), Lemma \ref{lem-1}, Theorem
\Ref{Green-thm} and Lebesgue's Dominated Convergence theorem, we see
that

\begin{eqnarray*}
M_{p}^{p}(r,D_{f})&=&\|D_{f}(0)\|^{p}+\frac{1}{2}\int_{\mathbb{D}(r)}\Delta\big(\|D_{f}(z)\|^{p}\big)\log\frac{r}{|z|}dA(z)\\
&=&\|D_{f}(0)\|^{p}\\
&&+\frac{1}{2}\int_{\mathbb{D}(r)}\Bigg\{p\|D_{f}(z)\|^{p-2}{\rm
Re}\Big[\overline{f_{z}(z)}(\Delta
f(z))_{z}+\overline{f_{\overline{z}}(z)}(\Delta
f(z))_{\overline{z}}\Big]\\
&&+p(p-2)\|D_{f}(z)\|^{p-4}\bigg|
\Big(f_{zz}(z)\overline{f_{z}(z)}\\
&&+\overline{f_{z\overline{z}}(z)}f_{z}(z)+
f_{\overline{z}z}(z)\overline{f_{\overline{z}}(z)}+
\overline{f_{\overline{z}\overline{z}}(z)}f_{\overline{z}}(z)\Big)\bigg|^{2}\\
&&+2p\|D_{f}(z)\|^{p-2}\big(D^{\ast}_{f}(z)\big)^{2}\Bigg\}\log\frac{r}{|z|}dA(z)\\
&\leq&\|D_{f}(0)\|^{p}+\frac{p}{2}\int_{\mathbb{D}(r)}\Big[(\eta+\tau)\|D_{f}(z)\|^{p}\\
&&+2(2p-3)\|D_{f}(z)\|^{p-2}\big(D^{\ast}_{f}(z)\big)^{2}
\Big]\log\frac{r}{|z|}dA(z)\\
&=&\|D_{f}(0)\|^{p}+p(\eta+\tau)\int_{0}^{r}M_{p}^{p}(\rho,
D_{f})\rho\log\frac{r}{\rho}d\rho+\\
&&2p(2p-3)\int_{0}^{r}\rho\log\frac{r}{\rho}\left(\frac{1}{2\pi}\int_{0}^{2\pi}\|D_{f}(\rho
e^{i\theta})\|^{p-2}
\big(D^{\ast}_{f}(\rho e^{i\theta})\big)^{2}d\theta\right) d\rho\\
&\leq&\|D_{f}(0)\|^{p}+p(\eta+\tau)\int_{0}^{r}M_{p}^{p}(\rho,
D_{f})\rho\log\frac{r}{\rho}d\rho\\
&&+2p(2p-3)\int_{0}^{r}M_{p}^{2}(\rho,D^{\ast}_{f})M_{p}^{p-2}(\rho,D_{f})\rho\log\frac{r}{\rho}\,d\rho\\
&\leq&\|D_{f}(0)\|^{p}+p(\eta+\tau)M_{p}^{p}(r,
D_{f})\int_{0}^{r}\rho\log\frac{r}{\rho}d\rho\\
&&+2p(2p-3)\int_{0}^{r}M_{p}^{2}(\rho,D^{\ast}_{f})M_{p}^{p-2}(\rho,D_{f})\rho\log\frac{r}{\rho}\,d\rho,
\end{eqnarray*}
which, together with (\ref{eq-CR-36}), 
gives that

\begin{eqnarray*}
&&\left[1-p(\eta+\tau)\int_{0}^{r}\rho\log\frac{r}{\rho}d\rho\right]M_{p}^{2}(r,D_{f})\\
&=&\left[1-\frac{pr^{2}(\eta+\tau)}{4}\right]M_{p}^{2}(r,D_{f})\\
&\leq&\|D_{f}(0)\|^{2}+2p(2p-3)\int_{0}^{r}\rho\log\frac{r}{\rho}M_{p}^{2}(\rho,D^{\ast}_{f})d\rho\\
&\leq&\|D_{f}(0)\|^{2}+2p(2p-3)\int_{0}^{r}(r-\rho)M_{p}^{2}(\rho,D^{\ast}_{f})d\rho\\
&=&\|D_{f}(0)\|^{2}+2p(2p-3)r^{2}\int_{0}^{1}(1-t)M_{p}^{2}\big(rt,D^{\ast}_{f}\big)dt\\
&\leq&\|D_{f}(0)\|^{2}+2p(2p-3)r^{2}(M^{\ast})^{2}\int_{0}^{1}\left[\omega\Big(\frac{1}{1-rt}\Big)\right]^{2}(1-t)dt\\
&\leq&\|D_{f}(0)\|^{2}+2p(2p-3)r^{2}(M^{\ast})^{2}\int_{0}^{1}\left[\omega\Big(\frac{1}{1-rt}\Big)\right]^{2}(1-rt)dt\\
&\leq&\|D_{f}(0)\|^{2}+2p(2p-3)r^{2}(M^{\ast})^{2}\omega(1)\int_{0}^{1}\omega\Big(\frac{1}{1-rt}\Big)dt,
\end{eqnarray*} where
 $dA$ denotes the
normalized area measure in $\mathbb{D}$. Then
\be\label{eq-CR-38}M_{p}^{2}(r,D_{f})\leq M^{\ast\ast}_{2}
\left[\|D_{f}(0)\|^{2}+M^{\ast}_{1}\int_{0}^{1}\omega\Big(\frac{1}{1-rt}\Big)dt\right],
\ee where $M^{\ast\ast}_{2}=1/\left[1-p(\eta+\tau)/4\right].$

\noindent ${\rm \mathbf{Step~ 4.}}$ By  (\ref{eq-CR37}),
(\ref{eq-CR36}),  Lemmas \ref{lem-4} and \ref{lem-5}, we obtain

\begin{eqnarray*}
M_{p}^{p}(r, f)&=&|f(0)|^{p}+\frac{1}{2}\int_{\mathbb{D}(r)}\Delta(|f(z)|^{p})\log\frac{r}{|z|}\,dV_{N}(z)\\
&\leq&|f(0)|^{p}+2p(p-1)\int_{0}^{r}\left(\frac{1}{2\pi}\int_{0}^{2\pi}|f(\rho
e^{i\theta})|^{p-2}
|D_{f}(\rho e^{i\theta})|^{2}d\theta\right)\rho\log\frac{r}{\rho}d\rho\\
&&+p(\eta+\tau)\int_{0}^{r}\left(\frac{1}{2\pi}\int_{0}^{2\pi}|f(\rho e^{i\theta})|^{p}\,d\theta\right)\rho\log\frac{r}{\rho}d\rho\\
&\leq&|f(0)|^{p}+2p(p-1)\int_{0}^{r}\rho\log\frac{r}{\rho}M_{p}^{2}(\rho,D_{f}) M_{p}^{p-2}(\rho, f)\,d\rho\\
&&+p(\eta+\tau)M_{p}^{p}(r,f)\int_{0}^{r}\rho\log\frac{r}{\rho}\,d\rho\\
&\leq&|f(0)|^{p}+2p(p-1)\int_{0}^{r}\rho\log\frac{r}{\rho}M_{p}^{2}(\rho,D_{f}) M_{p}^{p-2}(\rho, f)\,d\rho\\
&&+\frac{pr^{2}(\eta+\tau)}{4}M_{p}^{p}(r,f).
\end{eqnarray*}
The above, (\ref{eq-CR-36}), (\ref{eq-CR-38}) and the
monotonicity of $M_p(r,f)$ on $r$, imply that
\begin{eqnarray*}
 \frac{M_{p}^{2}(r,
 f)}{M^{\ast\ast}_{2}}
 &\leq&\left[1-\frac{p(\eta+\tau) r^{2}}{4}\right]M_{p}^{2}(r,
 f)\\
 &\leq&|f(0)|^{2}+2p(p-1)\int_{0}^{r}M_{p}^{2}(\rho,D_{f})\rho\log\frac{r}{\rho}\,d\rho\\
&\leq&|f(0)|^{2}+2p(p-1)\int_{0}^{r}M_{p}^{2}(\rho,D_{f})(r-\rho)\,d\rho\\
&\leq&|f(0)|^{2}+2p(p-1)\int_{0}^{1}M_{p}^{2}(r\rho,D_{f})(1-\rho)\,d\rho\\
&\leq&|f(0)|^{2}+2p(p-1)M^{\ast\ast}_{2}\|D_{f}(0)\|^{2}\\
&&+2p(p-1)M^{\ast}_{1}M^{\ast\ast}_{2}\int_{0}^{1}\left[\int_{0}^{1}\omega\Big(\frac{1}{1-r\rho
t}\Big)(1-\rho)dt\right]d\rho\\
&\leq&|f(0)|^{2}+2p(p-1)M^{\ast\ast}_{2}\|D_{f}(0)\|^{2}\\
&&+2p(p-1)M^{\ast}_{1}M^{\ast\ast}_{2}\int_{0}^{1}\left[\int_{0}^{1}\omega\Big(\frac{1}{1-r\rho
t}\Big)(1-rt\rho)dt\right]d\rho\\
&\leq&|f(0)|^{2}+2p(p-1)M^{\ast\ast}_{2}\|D_{f}(0)\|^{2}+2p(p-1)M^{\ast}_{1}M^{\ast\ast}_{2}\omega(1)\\
&<&\infty.
\end{eqnarray*}
Hence $f\in \mathcal{H}_{g}^{p}(\mathbb{D})$. The proof of the
theorem is complete. \qed

\section{Dirichlet-type spaces, Bergman-type spaces and applications to PDEs}\label{csw-sec3}
\subsection*{Proof of Theorem \ref{thm-2}}
 We first prove the necessity. Since ${\rm Re}(f\overline{\Delta f})\geq0$, 
we observe that $\Delta\big(|f|^{p}\big)\geq0$ and
\be\label{eq-CR15}0\leq\int_{\mathbb{B}^{n}}(1-|z|^{2})^{\alpha}\Delta\big(|f(z)|^{p}\big)dV_{N}(z)<\infty.\ee

Let  $r\in(0,1)$. For $ \alpha\geq2$, it is not difficult to see
that
$$(r^{2}-|z|^{2})^{\alpha}|_{\partial\mathbb{B}^{n}(r)}=0~\mbox{and}~\frac{\partial}{\partial
\varepsilon}\big[(r^{2}-|z|^{2})^{\alpha}\big]|_{\partial\mathbb{B}^{n}(r)}=0,
$$ where $\partial/\partial\varepsilon$ denotes an outer normal
derivative. Then, by Green's theorem, we get

\beq\label{eq-CR13}
\nonumber\int_{\mathbb{B}^{n}(r)}(r^{2}-|z|^{2})^{\alpha}\Delta\big(|f(z)|^{p}\big)dV_{N}(z)
&=&\int_{\mathbb{B}^{n}(r)}|f(z)|^{p}\Delta\big[(r^{2}-|z|^{2})^{\alpha}\big]dV_{N}(z)\\
\nonumber&=&4\alpha\int_{\mathbb{B}^{n}(r)}|f(z)|^{p}(r^{2}-|z|^{2})^{\alpha-2}\\
&&\times\left[|z|^{2}(n+\alpha-1)-nr^{2}\right]dV_{N}(z), \eeq
which, together with (\ref{eq-CR15}), gives that

\begin{eqnarray*}
\infty&>&4\alpha\int_{\mathbb{B}^{n}(R_{1})}
|f(z)|^{p}(1-|z|^{2})^{\alpha-2}\left[n-|z|^{2}(n+\alpha-1)\right]dV_{N}(z)\\
&&+\int_{\mathbb{B}^{n}}(1-|z|^{2})^{\alpha}\Delta\big(|f(z)|^{p}\big)dV_{N}(z)\\
 &\geq&4\alpha\int_{\mathbb{B}^{n}(rR_{1})}
|f(z)|^{p}(1-|z|^{2})^{\alpha-2}\left[n-|z|^{2}(n+\alpha-1)\right]dV_{N}(z)\\
&&+\int_{\mathbb{B}^{n}(r)}(1-|z|^{2})^{\alpha}\Delta\big(|f(z)|^{p}\big)dV_{N}(z)\\
&\geq&4\alpha\int_{\mathbb{B}^{n}(rR_{1})}
|f(z)|^{p}(r^{2}-|z|^{2})^{\alpha-2}\left[nr^{2}-|z|^{2}(n+\alpha-1)\right]dV_{N}(z)\\
&&+\int_{\mathbb{B}^{n}(r)}(r^{2}-|z|^{2})^{\alpha}\Delta\big(|f(z)|^{p}\big)dV_{N}(z)\\
&=&4\alpha\int_{\mathbb{B}^{n}(r)\setminus\mathbb{B}^{n}(rR_{1})}
|f(z)|^{p}(r^{2}-|z|^{2})^{\alpha-2}\left[|z|^{2}(n+\alpha-1)-nr^{2}\right]dV_{N}(z)\\
&\geq&2r^{2}\alpha\int_{\mathbb{B}^{n}(r)\setminus\mathbb{B}^{n}(rR_{2})}
|f(z)|^{p}(r^{2}-|z|^{2})^{\alpha-2}dV_{N}(z),
\end{eqnarray*} where
$R_{1}=\sqrt{\frac{n}{n+\alpha-1}}$ and
$R_{2}=\sqrt{\frac{n+\frac{1}{2}}{n+\alpha-1}}$.

 For
$R_{2}<r<1$,  we conclude that \beq\label{eq-CR16}
\nonumber\infty&>&2r^{2}\alpha\int_{\mathbb{B}^{n}(r)\setminus\mathbb{B}^{n}(rR_{2})}
|f(z)|^{p}(r^{2}-|z|^{2})^{\alpha-2}dV_{N}(z)\\
&\geq&2\alpha R_{2}^{2}U(r), \eeq
 where  $$U(r)=\int_{\mathbb{B}^{n}(r)\setminus\mathbb{B}^{n}(R_{2})}
|f(z)|^{p}(r^{2}-|z|^{2})^{\alpha-2}dV_{N}(z).$$ Then, for
$R_{2}<r<1$, $U(r)$ is increasing and bounded, from which we conclude that
$$\lim_{r\rightarrow1-}U(r)$$ exists. Hence for $p\geq2$,
 $f\in b_{\alpha-2,p}(\mathbb{B}^{n})$.

Next we prove the sufficiency. For $\alpha\geq2$, by
(\ref{eq-CR13}), we have

\beq \label{eq-CR17} &&4\alpha\int_{\mathbb{B}^{n}(rR_{1})}
|f(z)|^{p}(r^{2}-|z|^{2})^{\alpha-2}\left[nr^{2}-|z|^{2}(n+\alpha-1)\right]dV_{N}(z)\\
\nonumber
&&+\int_{\mathbb{B}^{n}(r)}(r^{2}-|z|^{2})^{\alpha}\Delta\big(|f(z)|^{p}\big)dV_{N}(z)\\
\nonumber
&=&4\alpha\int_{\mathbb{B}^{n}(r)\setminus\mathbb{B}^{n}(rR_{1})}
|f(z)|^{p}(r^{2}-|z|^{2})^{\alpha-2}\left[|z|^{2}(n+\alpha-1)-nr^{2}\right]dV_{N}(z)\\
\nonumber
&\leq&4\alpha(n+\alpha-1)\int_{\mathbb{B}^{n}(r)\setminus\mathbb{B}^{n}(rR_{1})}
|f(z)|^{p}(r^{2}-|z|^{2})^{\alpha-2}dV_{N}(z)\\ \nonumber
&\leq&4\alpha(n+\alpha-1)\int_{\mathbb{B}^{n}}|f(z)|^{p}(1-|z|^{2})^{\alpha-2}dV_{N}(z)\\
\nonumber &<&\infty. \eeq Since

\beq\label{eq-CR18} \nonumber\infty&>&\int_{\mathbb{B}^{n}(R_{1})}
|f(z)|^{p}(1-|z|^{2})^{\alpha-2}\left[n-|z|^{2}(n+\alpha-1)\right]dV_{N}(z)\\
&\geq&\int_{\mathbb{B}^{n}(rR_{1})}
|f(z)|^{p}(r^{2}-|z|^{2})^{\alpha-2}\left[nr^{2}-|z|^{2}(n+\alpha-1)\right]dV_{N}(z),
 \eeq which, together with (\ref{eq-CR17}) and $\Delta\big(|f|^{p}\big)\geq0$, implies that
$$\lim_{r\rightarrow1-}\int_{\mathbb{B}^{n}(r)}(r^{2}-|z|^{2})^{\alpha}\Delta\big(|f(z)|^{p}\big)dV_{N}(z)$$
does exist.
Therefore,$$\int_{\mathbb{B}^{n}}(1-|z|^{2})^{\alpha}\Delta\big(|f(z)|^{p}\big)dV_{N}(z)<\infty,$$
and thus the theorem is proved. \qed

The following result is well-known.
\begin{lem}\label{Lem-2}
Suppose that $a,b\in[0,\infty)$ and $q\in(0,\infty)$. Then
$$(a+b)^{q}\leq2^{\max\{q-1,0\}}(a^{q}+b^{q}).
$$
\end{lem}

\subsection*{Proof of Theorem \ref{thm-1}} By Lemma \ref{lem-1}, for $\rho\in[0,d(z))$, we get
\be\label{eq-CR1}\|D_{f}(z)\|^{\alpha}\leq\int_{\partial\mathbb{B}^{n}}\|D_{f}(z+\rho\zeta)\|^{\alpha}d\sigma(\zeta).\ee
Multiplying both sides of the inequality (\ref{eq-CR1}) by
$2n\rho^{2n-1}$ and integrating from $0$ to $d(z)/2$, we have
\begin{eqnarray*}
\frac{d(z)^{2n}\|D_{f}(z)\|^{\alpha}}{2^{2n}}&\leq&\int_{\partial\mathbb{B}^{n}}\int_{0}^{\frac{d(z)}{2}}2n\rho^{2n-1}
\|D_{f}(z+\rho\zeta)\|^{\alpha}d\rho d\sigma(\zeta)\\
&=&\int_{\mathbb{B}^{n}\big(z,\frac{d(z)}{2}
\big)}\|D_{f}(\xi)\|^{\alpha}dV_{N}(\xi)\\
&\leq&2^{\gamma}(d(z))^{-\gamma}\int_{\mathbb{B}^{n}\big(z,\frac{d(z)}{2}
\big)}(1-|\xi|)^{\gamma}\|D_{f}(\xi)\|^{\alpha}dV_{N}(\xi)\\
&\leq&\frac{2^{\gamma}\|f\|^{\alpha}_{\mathcal{D}_{\gamma,\alpha}}}{\big(d(z)\big)^{\gamma}},
\end{eqnarray*}
which implies that
\be\label{eq-CR2}\|D_{f}(z)\|\leq\frac{M_{1}}{(d(z))^{q+1}}, \ee
where $M_{1}=2^{1+q}\|f\|_{\mathcal{D}_{\gamma,\alpha}}$ and
$q=\frac{\gamma+2n}{\alpha}-1.$
By (\ref{eq-CR2}), 
 we know that

\beq\label{eq-CR3} \nonumber
|f(z)|&\leq&|f(0)|+\left|\int_{[0,z]}df(\varsigma)\right|\\
\nonumber&\leq&|f(0)|+\sqrt{2}\int_{[0,z]}\|D_{f}(\varsigma)\|\|d\varsigma\|\\
&\leq&|f(0)|+\frac{M_{2}}{\big(d(z)\big)^{q}}, \eeq where
$M_{2}=M_{1}\sqrt{2}/q$ and $[0,z]$ denotes the line segment from $0$ to
$z$.

By (\ref{eq-CR3}) and Lemma \ref{Lem-2}, we see that for
$z\in\mathbb{B}^{n}$,

\beq\label{eq-CR4}
|f(z)|^{p-2}&\leq&\left[|f(0)|+\frac{M_{2}}{\big(d(z)\big)^{q}}\right]^{p-2}
\leq
2^{p-2}\left[|f(0)|^{p-2}+\frac{M_{2}^{p-2}}{\big(d(z)\big)^{q(p-2)}}\right],
\eeq

\beq\label{eq-CR5}
|f(z)|^{p-1}&\leq&\left[|f(0)|+\frac{M_{2}}{\big(d(z)\big)^{q}}\right]^{p-1}
\leq
2^{p-1}\left[|f(0)|^{p-1}+\frac{M_{2}^{p-1}}{\big(d(z)\big)^{q(p-1)}}\right]
\eeq and \beq\label{eq-CR6}
|f(z)|^{p}&\leq&\left[|f(0)|+\frac{M_{2}}{\big(d(z)\big)^{q}}\right]^{p}
\leq
2^{p}\left[|f(0)|^{p}+\frac{M_{2}^{p}}{\big(d(z)\big)^{qp}}\right].
\eeq

\noindent ${\rm \mathbf{Case~ 1.}}$ Let $p\in[4,\infty)$.

By direct calculations, we get

\beq\label{eq-CR7}
\nonumber\Delta\big(|f|^{p}\big)&=&p(p-2)|f|^{p-4}\sum_{k=1}^{n}|f_{z_{k}}\overline{f}+\overline{f_{\overline{z}_{k}}}f|^{2}+
2p|f|^{p-2}\|D_{f}\|^{2}+p|f|^{p-2}\mbox{Re}(f\overline{\Delta f})\\
\nonumber
&\leq&p(p-2)|f|^{p-4}\sum_{k=1}^{n}|f_{z_{k}}\overline{f}+\overline{f_{\overline{z}_{k}}}f|^{2}+
2p|f|^{p-2}\|D_{f}\|^{2}+p|f|^{p-1}|\Delta f|\\
&\leq&2p(p-1)|f|^{p-2}\|D_{f}\|^{2}+pa|f|^{p-1}\|D_{f}\|+pb|f|^{p}+pc|f|^{p-1}.
\eeq It follows from (\ref{eq-CR4}), (\ref{eq-CR5}), (\ref{eq-CR6})
and (\ref{eq-CR7}) that

\beq\label{eq-CR8}
\nonumber\big(d(z)\big)^{pq}\Delta\big(|f|^{p}\big)
&\leq&2p(p-1)\big(d(z)\big)^{pq}|f|^{p-2}\|D_{f}\|^{2}\\
\nonumber &&+pa\big(d(z)\big)^{pq}|f|^{p-1}\|D_{f}\|+
pb\big(d(z)\big)^{pq}|f|^{p}\\
\nonumber
&&+pc\big(d(z)\big)^{pq}|f|^{p-1}\\
\nonumber&=&2p(p-1)\big(d(z)\big)^{pq-\frac{2\gamma}{\alpha}}
|f|^{p-2}\|D_{f}\|^{2}\big(d(z)\big)^{\frac{2\gamma}{\alpha}}\\
\nonumber &&+pa\big(d(z)\big)^{pq-\frac{\gamma}{\alpha}}|f|^{p-1}
\|D_{f}\|\big(d(z)\big)^{\frac{\gamma}{\alpha}}\\
\nonumber&&+
pb\big(d(z)\big)^{pq}|f|^{p}+pc\big(d(z)\big)^{pq}|f|^{p-1}\\
&\leq&M_{3}\|D_{f}\|^{2}\big(d(z)\big)^{\frac{2\gamma}{\alpha}}+M_{4}\|D_{f}\|\big(d(z)\big)^{\frac{\gamma}{\alpha}}+M_{5},
\eeq where
$$M_{3}=2^{p-1}p(p-1)\left(|f(0)|^{p-2}+M_{2}^{p-2}\right),$$ $$M_{4}=p2^{p-1}\left(|f(0)|^{p-1}+M_{2}^{p-1}\right)
\sup_{z\in\mathbb{B}^{n}}a(z)$$ and
$$M_{5}=p2^{p}\left(|f(0)|^{p}+M_{2}^{p}\right)\sup_{z\in\mathbb{B}^{n}}b(z)+pM2^{p-1}\left(|f(0)|^{p-1}+M_{2}^{p-1}\right).$$
By H\"older's inequality, we obtain \beq\label{eq-CR9}
\int_{\mathbb{B}^{n}}\big(d(z)\big)^{\frac{2\gamma}{\alpha}}\|D_{f}(z)\|^{2}dV_{N}(z)&\leq&
\left(\int_{\mathbb{B}^{n}}\big(d(z)\big)^{\gamma}\|D_{f}(z)\|^{\alpha}dV_{N}(z)\right)^{\frac{2}{\alpha}}\\
\nonumber
&&\times\left(\int_{\mathbb{B}^{n}}dV_{N}(z)\right)^{1-\frac{2}{\alpha}}\\
\nonumber &\leq&\|f\|^{2}_{\mathcal{D}_{\gamma,\alpha}}, \eeq which
gives \beq\label{eq-CR10}
\int_{\mathbb{B}^{n}}\big(d(z)\big)^{\frac{\gamma}{\alpha}}\|D_{f}(z)\|dV_{N}(z)&\leq&
\left(\int_{\mathbb{B}^{n}}\big(d(z)\big)^{\frac{2\gamma}{\alpha}}\|D_{f}(z)\|^{2}dV_{N}(z)\right)^{\frac{1}{2}}\\
\nonumber
&&\times\left(\int_{\mathbb{B}^{n}}dV_{N}(z)\right)^{\frac{1}{2}}\\
\nonumber &\leq&\|f\|_{\mathcal{D}_{\gamma,\alpha}}. \eeq It follows
from   (\ref{eq-CR8}), (\ref{eq-CR9}) and (\ref{eq-CR10}) that

\begin{eqnarray*}
\int_{\mathbb{B}^{n}}\big(d(z)\big)^{pq}\Delta\big(|f(z)|^{p}\big)dV_{N}(z)&=&\int_{\mathbb{B}^{n}}
\Big[M_{3}\|D_{f}\|^{2}\big(d(z)\big)^{\frac{2\gamma}{\alpha}}\\
&&+M_{4}\|D_{f}\|\big(d(z)\big)^{\frac{\gamma}{\alpha}}+M_{5}\Big]dV_{N}(z)\\
&\leq&M_{3}\|f\|^{2}_{\mathcal{D}_{\gamma,\alpha}}+M_{4}\|f\|_{\mathcal{D}_{\gamma,\alpha}}+M_{5}\\
&<&\infty.
\end{eqnarray*}

\noindent ${\rm \mathbf{Case~ 2.}}$ Let $p\in[2,4)$.

For  $p\in[2,4)$, $m\in\{1,2,\ldots\}$ and $r\in(0,1)$, let
$T_{m}^{p}=\left(|f|^{2}+\frac{1}{m}\right)^{\frac{p}{2}}$. Then, by
(\ref{eq-CR8}), (\ref{eq-CR9}), (\ref{eq-CR10}) and Lebesgue's
Dominated Convergence Theorem, we have

\begin{eqnarray*}
&&\lim_{r\rightarrow1-}\left\{\lim_{m\rightarrow\infty}\int_{\mathbb{B}^{n}(r)}\big(d(z)\big)^{pq}
\Delta\big(T_{m}^{p}(z)\big)dV_{N}(z)\right\}\\
&=&\lim_{r\rightarrow1-}\int_{\mathbb{B}^{n}(r)}\big(d(z)\big)^{pq}
\lim_{m\rightarrow\infty}\Delta\big(T_{m}^{p}(z)\big)dV_{N}(z)\\
&=&\lim_{r\rightarrow1-}\int_{\mathbb{B}^{n}(r)}\big(d(z)\big)^{pq}
\bigg[p(p-2)|f(z)|^{p-4}\sum_{k=1}^{n}\big|f_{z_{k}}(z)\overline{f(z)}\\
&&+\overline{f_{\overline{z}_{k}}(z)}f(z)\big|^{2}+
2p|f(z)|^{p-2}\|D_{f}(z)\|^{2}\\
&&+p|f(z)|^{p-2}\mbox{Re}(f(z)\overline{\Delta
f(z)})\bigg]dV_{N}(z)\\
&\leq&\int_{\mathbb{B}^{n}}\Big(M_{3}\|D_{f}\|^{2}\big(d(z)\big)^{\frac{2\gamma}{\alpha}}+
M_{4}\|D_{f}\|\big(d(z)\big)^{\frac{\gamma}{\alpha}}+M_{5}\Big)dV_{N}(z)\\
&<&\infty.
\end{eqnarray*}
This concludes the proof of the theorem. \qed

\subsection*{Proof of Theorem \ref{thm-4}} \noindent ${\rm \mathbf{Case~ 1.}}$ Let $n\geq2$.

 Without loss of
generality, we may assume that
$$\inf_{z\in\mathbb{B}^{n}}a_{1}(z)>0~\mbox{ and}~
\inf_{z\in\mathbb{B}^{n}}b_{1}(z)>0.$$ Let $r_{0}\in(0,1)$ be a
constant. Then, by Lemma \ref{lem-5}, for $0<r_{0}\leq r<1$, we have

\beq\label{eq-CR31} M_{p}^{p}(r,f)
&=&|f(0)|^{p}+\int_{\mathbb{B}^{n}(r)}\Delta\big(|f(z)|^{p}\big)G_{2n}(z,r)dV_{N}(z)\\
\nonumber
&\geq&|f(0)|^{p}+\int_{\mathbb{B}^{n}(r)}\Big(a_{1}(z)\|D_{f}(z)\|^{t_{1}}+b_{1}(z)|f(z)|^{t_{2}}\\
\nonumber &&+c_{1}(z)\Big)G_{2n}(z,r)dV_{N}(z)\\ \nonumber \eeq \beq
\nonumber &\geq&
|f(0)|^{p}+\inf_{z\in\mathbb{B}^{n}}a_{1}(z)\int_{\mathbb{B}^{n}(r)}\|D_{f}(z)\|^{t_{1}}G_{2n}(z,r)dV_{N}(z)\\
\nonumber
&&+\inf_{z\in\mathbb{B}^{n}}b_{1}(z)\int_{\mathbb{B}^{n}(r)}|f(z)|^{t_{2}}G_{2n}(z,r)dV_{N}(z)\\
\nonumber
&&+\inf_{z\in\mathbb{B}^{n}}c_{1}(z)\int_{\mathbb{B}^{n}(r)}G_{2n}(z,r)dV_{N}(z)\\
\nonumber
\nonumber&=&|f(0)|^{p}+\inf_{z\in\mathbb{B}^{n}}a_{1}(z)\int_{\mathbb{B}^{n}(r_{0})}\|D_{f}(z)\|^{t_{1}}G_{2n}(z,r)dV_{N}(z)\\
\nonumber &&+
\inf_{z\in\mathbb{B}^{n}}a_{1}(z)\int_{\mathbb{B}^{n}(r)\setminus\mathbb{B}^{n}(r_{0})}\|D_{f}(z)\|^{t_{1}}G_{2n}(z,r)dV_{N}(z)\\
 \nonumber
&&+\inf_{z\in\mathbb{B}^{n}}b_{1}(z)\int_{\mathbb{B}^{n}(r_{0})}|f(z)|^{t_{2}}G_{2n}(z,r)dV_{N}(z)\\
 \nonumber
&&+\inf_{z\in\mathbb{B}^{n}}b_{1}(z)\int_{\mathbb{B}^{n}(r)\setminus\mathbb{B}^{n}(r_{0})}|f(z)|^{t_{2}}G_{2n}(z,r)dV_{N}(z)\\
\nonumber
&&+\inf_{z\in\mathbb{B}^{n}}c_{1}(z)\int_{\mathbb{B}^{n}(r)}G_{2n}(z,r)dV_{N}(z).
\eeq

It is easy to see  that, for all $r\in(0,1)$,
\be\label{eq-CR32}0<\int_{\mathbb{B}^{n}(r)}G_{2n}(z,r)dV_{N}(z)<\infty.\ee

Since \beq \nonumber
\int_{\mathbb{B}^{n}(r_{0})}\|D_{f}(z)\|^{t_{1}}G_{2n}(z,r)dV_{N}(z)&\leq&
\int_{\mathbb{B}^{n}(r_{0})}\|D_{f}(z)\|^{t_{1}}G_{2n}(z,1)dV_{N}(z)\\
\nonumber &<&\infty \eeq and

\beq
\nonumber\infty&>&\int_{\mathbb{B}^{n}(r)\setminus\mathbb{B}^{n}(r_{0})}\|D_{f}(z)\|^{t_{1}}G_{2n}(z,r)\big)dV_{N}(z)\\
\nonumber
&=&\frac{1}{4n(n-1)}\int_{\mathbb{B}^{n}(r)\setminus\mathbb{B}^{n}(r_{0})}\frac{(r-|z|)\big(\sum_{k=0}^{2n-3}r^{2n-3-k}|z|^{k}\big)}{|z|^{2n-2}r^{2n-2}}
\|D_{f}(z)\|^{t_{1}}dV_{N}(z)\\ \nonumber &\geq&\delta(r), \eeq
which, together with $f\in\mathcal{H}_{g}^{p}(\mathbb{B}^{n})$,
(\ref{eq-CR31}),  (\ref{eq-CR32}) and  the monotonicity of
$\delta(r)$, yield that the limit
\be\label{eq-CR29}\lim_{r\rightarrow1-}\int_{\mathbb{B}^{n}(r)\setminus\mathbb{B}^{n}(r_{0})}(r-|z|)\|D_{f}(z)\|^{t_{1}}dV_{N}(z)\ee
 exists, where
$$\delta(r)=\frac{r_{0}^{2n-3}}{2n}\int_{\mathbb{B}^{n}(r)\setminus\mathbb{B}^{n}(r_{0})}(r-|z|)\|D_{f}(z)\|^{t_{1}}dV_{N}(z).$$
Then $f\in\mathcal{D}_{1,t_{1}}(\mathbb{B}^{n})$.

By using a similar argument as in the proof of (\ref{eq-CR29}), we
see that $ f\in b_{1,t_{2}}(\mathbb{B}^{n})$.

\noindent ${\rm \mathbf{Case~ 2.}}$ Let $n=1$.

In this case, the proof is similar to the proof of  the case 2 in
Theorem \ref{thm-3}. Therefore, proof of the theorem is complete. \qed

\subsection*{Proof of Corollary \ref{thm-4.0}}
Without loss of generality, we assume that
$\prod_{k=1}^{n}\lambda_{k}\neq0.$

 \noindent ${\rm \mathbf{Case~ 1.}}$ Let
$p\in[4,\infty)$.

By computations, for $k\in\{1,\ldots,n\}$, we have
\begin{eqnarray*}
(|f|^{p})_{z_{k}\overline{z}_{k}}&=&[(f^{\frac{p}{2}}\overline{f}^{\frac{p}{2}})_{\overline{z}_{k}}]_{z_{k}}\\
&=&\frac{p}{2}\left(f^{\frac{p}{2}-1}f_{\overline{z}_{k}}\overline{f}^{\frac{p}{2}}+
f^{\frac{p}{2}}\overline{f}^{\frac{p}{2}-1}\overline{f_{z_{k}}}\right)_{z_{k}}\\
&=&\frac{p}{2}\left(\lambda_{k}
f^{\frac{p}{2}+\frac{\alpha}{2}-1}\overline{f}^{\frac{p}{2}+\frac{\alpha}{2}}+
f^{\frac{p}{2}}\overline{f}^{\frac{p}{2}-1}\overline{f_{z_{k}}}\right)_{z_{k}}\\
&=&\frac{p}{2}\Big[\lambda_{k}\big(\frac{p}{2}+\frac{\alpha}{2}-1\big)f^{\frac{p}{2}+
\frac{\alpha}{2}-2}\overline{f}^{\frac{p}{2}+\frac{\alpha}{2}}f_{z_{k}}+
\lambda_{k}\big(\frac{p}{2}+\frac{\alpha}{2}\big)|f|^{p+\alpha-2}\overline{f_{\overline{z}_{k}}}\\
&&+\frac{p}{2}|f|^{p-2}|f_{z_{k}}|^{2}+
\big(\frac{p}{2}-1\big)\overline{f_{z_{k}}}f^{\frac{p}{2}}\overline{f}^{\frac{p}{2}-2}\overline{f_{\overline{z}_{k}}}
+f^{\frac{p}{2}}\overline{f}^{\frac{p}{2}-1}\overline{f_{z_{k}\overline{z}_{k}}}\Big]\\
&=&\frac{p}{2}\Big[\lambda_{k}\big(\frac{p}{2}+\frac{\alpha}{2}-1\big)|f|^{p+\alpha-4}f_{z_{k}}\overline{f}^{2}+
\lambda^{2}_{k}\big(\frac{p+\alpha}{2}\big)|f|^{p+\alpha-2}|f|^{\alpha}\\
&&+\frac{p}{2}|f|^{p-2}|f_{z_{k}}|^{2}
+\lambda_{k}\big(\frac{\alpha}{2}+\frac{p}{2}-1\big)|f|^{p+\alpha-4}\overline{f_{z_{k}}}f^{2}+\frac{\alpha\lambda_{k}^{2}}{2}|f|^{p+2\alpha-2}\Big]\\
&=&\frac{p}{2}\Big[\lambda_{k}\big(\frac{p}{2}+\frac{\alpha}{2}-1\big)|f|^{p+\alpha-4}f_{z_{k}}\overline{f}^{2}+
\lambda_{k}\big(\frac{p}{2}+\frac{\alpha}{2}-1\big)|f|^{p+\alpha-4}\overline{f_{z_{k}}}f^{2}\\
&&+\lambda_{k}^{2}\big(\alpha+\frac{p}{2}\big)|f|^{p+2\alpha-2}+\frac{p}{2}|f_{z_{k}}|^{2}|f|^{p-2}\Big]\\
&=&\frac{p}{2}\Big\{\mbox{Re}\big[\lambda_{k}(p+\alpha-2)|f|^{p+\alpha-4}f_{z_{k}}\overline{f}^{2}\big]+
\lambda_{k}^{2}\big(\alpha+\frac{p}{2}\big)|f|^{p+2\alpha-2}\\
&&+\frac{p}{2}|f_{z_{k}}|^{2}|f|^{p-2}\Big\},
\end{eqnarray*}
which implies
\begin{eqnarray*}
\Delta(|f|^{p})&=&4\sum_{k=1}^{n}(|f|^{p})_{z_{k}\overline{z}_{k}}\\
&=&2p\Big\{\sum_{k=1}^{n}\mbox{Re}\big[\lambda_{k}(p+\alpha-2)|f|^{p+\alpha-4}f_{z_{k}}\overline{f}^{2}\big]\\
&&+
\sum_{k=1}^{n}\lambda_{k}^{2}\big(\alpha+\frac{p}{2}\big)|f|^{p+2\alpha-2}+\frac{p}{2}\sum_{k=1}^{n}|f_{z_{k}}|^{2}|f|^{p-2}\Big\}.
\end{eqnarray*}
Hence,
\begin{eqnarray*}
\Delta(|f|^{p})&-&[4p-(2-\alpha)^{2}]|f|^{p+2\alpha-2}\sum_{k=1}^{n}\lambda^{2}_{k}\\
\end{eqnarray*}
\begin{eqnarray*}
&=&\sum_{k=1}^{n}\Big\{\lambda_{k}^{2}\big[p^{2}+2p(\alpha-2)+(\alpha-2)^{2}\big]|f|^{p+2\alpha-2}
\\
&&+\mbox{Re}\big[2p\lambda_{k}(p+\alpha-2)|f|^{p+\alpha-4}f_{z_{k}}\overline{f}^{2}\big]
+p^{2}|f_{z_{k}}|^{2}|f|^{p-2}\Big\}\\
&\geq&\sum_{k=1}^{n}\Big\{\lambda_{k}^{2}(p+\alpha-2)^{2}|f|^{p+2\alpha-2}+p^{2}|f_{z_{k}}|^{2}|f|^{p-2}\\
&&-|2p\lambda_{k}(p+\alpha-2)||f|^{p+\alpha-2}|f_{z_{k}}|\Big\}\\
&=&|f|^{p-2}\sum_{k=1}^{n}\big(|\lambda_{k}(p+\alpha-2)||f|^{\alpha}-|p||f_{z_{k}}|\big)^{2}\\
&\geq&0,
\end{eqnarray*} 
which yields

\be\label{eq-CR-41} \Delta(|f|^{p})\geq
\big[4p-(2-\alpha)^{2}\big]|f|^{p+2\alpha-2}\sum_{k=1}^{n}\lambda_{k}^{2}.\ee

 \noindent ${\rm \mathbf{Case~ 2.}}$
 Let $p\in[2,4)$.

For $m\in\{1,2,\ldots\}$, let
$G_{m}^{p}=\left(|f|^{2}+\frac{1}{m}\right)^{\frac{p}{2}}.$  Then by
 Lebesgue's Dominated Convergence theorem and (\ref{eq-CR-41}), we have

\begin{eqnarray*}
\lim_{m\rightarrow\infty}\Delta\left(G_{m}^{p}\right)&=&\lim_{m\rightarrow\infty}4\sum_{k=1}^{n}(G_{m}^{p})_{z_{k}\overline{z}_{k}}\\
&=&2p\Big\{\sum_{k=1}^{n}\mbox{Re}\big[\lambda_{k}(p+\alpha-2)|f|^{p+\alpha-4}f_{z_{k}}\overline{f}^{2}\big]\\
&&+
\sum_{k=1}^{n}\lambda_{k}^{2}\big(\alpha+\frac{p}{2}\big)|f|^{p+2\alpha-2}+\frac{p}{2}\sum_{k=1}^{n}|f_{z_{k}}|^{2}|f|^{p-2}\Big\}\\
 &\geq&
\big[4p-(2-\alpha)^{2}\big]|f|^{p+2\alpha-2}\sum_{k=1}^{n}\lambda_{k}^{2}.
\end{eqnarray*}

Applying Theorem \ref{thm-4}, we conclude that $f\in
b_{1,\vartheta}(\mathbb{B}^{n})$, where $\vartheta=p+2\alpha-2.$ The
proof of the Corollary is complete. \qed

\bigskip

{\bf Acknowledgements:} This research was partly  supported by
National Natural Science Foundation of China (No. 11401184 and No.
11071063), the Construct Program of the Key Discipline in Hunan
Province, and the V\"ais\"al\"a Foundation of the Finnish Academy of
Science and Letters.

\end{document}